\documentclass[12pt,a4paper]{article}
\usepackage[tbtags]{amsmath}    
\usepackage{amsfonts,amssymb,amsthm,euscript,makeidx,pifont,boxedminipage,multicol,fullpage}
\usepackage{color}

\def\be{\begin{eqnarray}}
\def\ee{\end{eqnarray}}

\def\b*{\begin{eqnarray*}}
\def\e*{\end{eqnarray*}}

\newtheorem{Theorem}{Theorem}[section]
\newtheorem{Lemma}[Theorem]{Lemma}
\newtheorem{Proposition}[Theorem]{Proposition}

\newtheorem{Definition}[Theorem]{Definition}
\newtheorem{Remark}[Theorem]{Remark}
\newtheorem{Example}[Theorem]{Example}

\newtheorem{Assumption}[Theorem]{Assumption}

\newcommand{\rmi}{{\rm (i)$\>\>$}}
\newcommand{\rmii}{{\rm (ii)$\>\>$}}
\newcommand{\rmiii}{{\rm (iii)$\>\>$}}
\newcommand{\rmiv}{{\rm (iv)$\>\>$}}

\def \E{\mathbb{E}}
\def \F{\mathbb{F}}

\def \L{\mathbb{L}}
\def \P{\mathbb{P}}
\def \Q{\mathbb{Q}}
\def \R{\mathbb{R}}

\def \N{\mathbb{N}}

\def\Ac{{\cal A}}
\def\Bc{{\cal B}}

\def\Fc{{\cal F}}
\def\Gc{{\cal G}}
\def\Hc{{\cal H}}

\def\Pc{{\cal P}}
\def\Qc{{\cal Q}}

\def\Uc{{\cal U}}

\def\Pt{{\tilde P}}

\def \Om{\Omega}
\def \om{\omega}
\def \Omb{\overline{\Om}}
\def \omb{\bar{\om}}
\def \omh{\hat{\om}}

\def \0{\mathbf{0}}
\def \1{\mathbf{1}}
\def \x{\times}

\def \Fcb{\overline{{\cal F}}}
\def \Fbb{\overline{\F}}

\def \Pcb{\overline{\Pc}}
\def \Qcb{\overline{\Qc}}
\def \Pt{\tilde{\P}}

\def \Pb{\overline{\P}}

\def \Hcb{\overline{\Hc}}
\def \Hb{\overline{H}}

\def\Bf{\mathfrak{B}}
\def\Rb{\overline{\R}}

\def\Qb{\overline{\Q}}
\def\Hcb{\overline{\Hc}}

\def\NA2{\mathrm{NA}_2}
\def\NA{\mathrm{NA}}

\def\Omh{\widehat{\Om}}
\def\Pch{\widehat{\Pc}}
\def\Pct{\widetilde{\Pc}}
\def\Qch{\widehat{\Qc}}
\def\Fch{\widehat{\Fc}}
\def\Fbh{\widehat{\F}}
\def\Ph{\widehat{\P}}
\def\Qh{\widehat{\Q}}
\def\cP{{\mathcal P}}

\def\CPSt{\mathcal{S}}

\def\scap#1#2{#1\cdot #2}

\title{Super-replication with  proportional transaction cost 
	 under model uncertainty}

\author{Bruno Bouchard\footnote{Universit\'e Paris-Dauphine, PSL Research University, CNRS, CEREMADE, Paris.}  \thanks{bouchard@ceremade.dauphine.fr} 
\and Shuoqing Deng\addtocounter{footnote}{-1}\footnotemark[\value{footnote}]~\addtocounter{footnote}{1}\thanks{deng@ceremade.dauphine.fr}
\and Xiaolu Tan\addtocounter{footnote}{-2}\footnotemark[\value{footnote}]~\addtocounter{footnote}{2}\thanks{tan@ceremade.dauphine.fr} \footnote{ Xiaolu Tan gratefully acknowledges the financial support of the ERC 321111 Rofirm, the ANR Isotace, and the Chairs Financial Risks (Risk Foundation, sponsored by Soci\'et\'e G\'en\'erale) and Finance and Sustainable Development (IEF sponsored by EDF and CA).}
}

\date{\today}

\begin{document}
\bibliographystyle{plain}

\maketitle

\abstract{
	We consider a discrete time financial market with proportional transaction cost under model uncertainty,
	and study a super-replication problem.
	We recover the duality results that are well known in the classical dominated context.
	Our key argument consists in using a randomization technique together with the minimax theorem  
	to convert the initial problem to a frictionless problem set on an enlarged space. This  
	allows  us to appeal to the techniques and results  of Bouchard and Nutz \cite{BouchardNutz.13}
	  to obtain the duality result.
}

\vspace{2mm}


\vspace{2mm}





\section{Introduction}

Discrete time  financial markets have been widely studied, and are now well understood. In the frictionless setting, the Fundamental Theorem of Asset Pricing and the traditional dual formulation of the set of contingent claims that can be super-hedged  from a zero initial endowment are proved by first showing that this latter set is closed in probability and by then using a Hahn-Banach separation argument in $\L^{1}$, as in the celebrated Kreps--Yan theorem~\cite{Yan.80}, see e.g.~\cite{DelbaenSchachermayer.06,FollmerSchied.04,bouchard2016fundamentals}.  In the presence of proportional transaction costs, serveral notions of no-arbitrage properties can be considered, but they  all aim at obtaining a similar closure property, so that Hahn-Banach separation arguments can still be  applied, see \cite{KabanovRasonyiStricker.02,Rasonyi.09,Schachermayer.04} and \cite{KabanovSafarian.09} for a survey monograph. 

In the context of model uncertainty, the market is defined with respect to a family $\cP$ of (typically singular) probability measures. Closure properties can still be proved, in the quasi-sure sense, but no satisfactory $\L^{p}$-type duality argument can be used, because of the lack of a reference (dominating) probability measure. 
 In the frictionless context, \cite{BouchardNutz.13} suggested to use a one period argument {\sl \`a la} Dalang--Morton--Willinger \cite{DalangMortonWillinger.90}, and then to appeal to measurable selection techniques to paste the periods together. This is unfortunatly not possible (in general) when proportional transaction costs are present: local no-arbitrage is not equivalent to global no-arbitrage. Still, a quasi-sure versions of the classical {\sl weak} and {\sl strict} no-arbitrage conditions could be characterized in \cite{BayraktarZhang.13}. This requires the use of an  intricate forward-backward construction.  A quite similar construction was later used in   Burzoni \cite{Burzoni16} for a version of the no-model independent arbitrage condition based on the robust no-arbitrage property of Schachermayer \cite{Schachermayer.04} .   
 
 Bouchard and Nutz \cite{bouchard2016consistent} proposed to follow a simpler route and to use a quasi-sure version of the only no-abitrage condition for which local no-arbitrage and global no-arbitrage are equivalent. This notion was first suggested by \cite{Rasonyi.09} under the name of {\sl no-arbitrage of the second kind}, or  {\sl no-sure gains in liquidation value}. Because of the equivalence between absence of local and absence of global arbitrage, they could use the same one period based arguments as in \cite{BouchardNutz.13} to provide a quasi-sure version of the Fundamental Theorem of Asset Pricing in this context.

Unfortunately, they were no able to come up with an equaly easy proof of the super-hedging duality that seems to require a global argument. The difficulty comes from the fact that a portfolio is described by a vector valued process. At time $t$, one needs to define a vector position allowing to super-hedge the required time $t+1$-position. In frictionless markets, the time $t+1$-position reduces to a scalar valued random variable, the time $t+1$-value of the super-hedging price. Its computation can be done backward, by iterating on the time periods.  In models with proportional transaction costs, there is an infinity of possible positions at time $t+1$, that are all minimal in the sense that none of them is dominated by another one, and that are enough to build up a super-hedging strategy.  A backward induction does not tell which one is consistent with the global super-hedging strategy.  

	Burzoni \cite{Burzoni16} was able to solve this issue by constructing a fictitious price system  in which  the frictionless superhedging price  is the same as in the original market.  In his context, the super-hedging has to hold pointwise on the so-called  {\sl efficient support of the family of consistent price systems}. Also the market is multivariate,  all transactions goes through the cash account (no direct exchanges between assets).   It complements the work of Dolinsky and Soner \cite{DolinskySoner.13} in which pointwise super-hedging is considered, and a static position on an arbitrary European  option can be initially taken. 

	\vspace{2mm}

	In this paper, we suggest to use a very simple randomization argument to tackle this problem in a general multivariate setting, under the quasi-sure no-arbitrage of second type condition of \cite{bouchard2016consistent}.
	In particular, we allow for direct exchanges between the risky assets, which is not possible in  \cite{Burzoni16,DolinskySoner.13}.
	Technically, we add additional randomness to our initial probability space to construct a {\sl fictitious} price process $X$ that is consistent with the original bid-ask bounds. This  additional randomness controls the positions of prices within these bounds. We then consider the problem of quasi-sure super-hedging in the fictionless market with price process $X$ and show that it matches with the super-hedging price in the original market with proportional transaction costs. This essentially follows from a  minimax argument in the one period setting, that can be iterated in a backward way. Then, it suffices to apply the super-hedging duality of Bouchard and Nutz \cite{BouchardNutz.13}, and to project back all involved quantities on the original probability space.
 
Note that this randomization/enlargement technique is 
	in fact in the same spirit of the {\sl controlled fictitious market} approach of \cite{bouchard2000explicit,cvitanic1999closed}, used in a dominated Markovian continuous time setting.

	\vspace{1mm}

	The rest of the paper is organized as follows. 
	In Section \ref{sec:main}, we  first describe our discrete time market with proportional transaction costs, and   introduce our randomization approach.
	We then specialize  to the probabilistic setting  suggested by Bouchard and Nutz \cite{bouchard2016consistent} and link   their  quasi-sure no-arbitrage condition of second kind to a quasi-sure no-arbitrage condition set on our randomized frictionless market.
	In Section  \ref{sec:pricing_hedging}, we   consider the super-replication problem and prove the duality, by using our randomization technique.

\vspace{2mm}

	\noindent {\bf Notations.} Given a measurable space $(\Om, \Fc)$, we denote by $\Bf(\Om, \Fc)$ the set of all probability measures on $(\Om, \Fc)$.
	If $\Om$ is a topological space, $\Bc(\Om)$ denotes its Borel $\sigma$-field and we abbreviate the notation 
	$\Bf(\Om) := \Bf(\Om, \Bc(\Om))$.
	If $\Om$ is a Polish space, a subset $A \subseteq \Om$ is analytic if it is the image of a Borel subset of another Polish space under a Borel measurable mapping.
	A function $f: \Om \to \Rb := [-\infty, \infty]$ is upper semianalytic if $\{\om \in \Om ~: f(\om) > a\}$ is analytic for all $a\in \R$.
	Given a probability measure $\P \in \Bf(\Om)$ and a measurable function $f: \Om \to \Rb$, we define the expectation
	$$
		\E^\P[ f] ~:=~ \E^{\P}[f^+] - \E^{\P}[f^-], ~~~\mbox{with the convention}~ \infty - \infty = -\infty.
	$$
	For a family $\Pc \subseteq \Bf(\Om)$ of probability measures, a subset $A \subset \Om$ is called $\Pc$-polar if $A \subset A'$ for some universally measurable set $A'$ satisfying $\P[A'] = 0$ for all $\P \in \Pc$,
	and a property is said to hold $\Pc$-quasi surely or $\Pc$-q.s if it holds true outside a $\Pc$-polar set.
	For $\Q \in \Bf(\Om)$, we write $\Q \lll \Pc$ if there exists $\P' \in \Pc$ such that $\Q \ll \P'$. Given a sigma algebra $\Gc$, we denote by $L^{0}(\Gc)$ the collection of $\R^{d}$-valued random variable that are $\Gc$-measurable, $d$ being given by the context.   If we are given a measurable random set $A$ and a family of probability measures $\Pc$, we denote by $L^{0}_{\Pc}(\Gc,A)$ the collection of $\Gc$-measurable random variables  taking values in $A$ $\Pc$-q.s.

\section{A randomization approach for market with proportional transaction cost} \label{sec:main}

	We first introduce an abstract discrete-time market with proportional transaction cost,
	and  show how to reduce to a fictitious market without transaction cost by using a randomization technique.
	Then, we specialize to the setting of Bouchard and Nutz \cite{bouchard2016consistent} and discuss in particular  how  their quasi-sure version of the second kind no-arbitrage condition can be related to a no-arbitrage condition set on our enlarged fictitious market.

\subsection{The financial market with proportional transaction cost}

	Let $(\Om, \Fc)$ be a measurable space, equipped with two filtrations $\F^{0} = (\Fc^{0}_t)_{t = 0, 1, \cdots, T}$ $\subset$   $\F = (\Fc_t)_{t = 0, 1, \cdots, T}$ for some $T \in \N$ (later on the first one will be the raw filtration, while the second one will be its universal completion).
	We fix a family of probability measures $\Pc$ on $(\Om, \Fc)$, which represents the model uncertainty.
	In particular, when $\Pc$ is singular, it reduces to the classical dominated market model framework.

	Following \cite{bouchard2016consistent}, we specify our financial market with proportional transaction cost in terms of random cones.
	Let $d \ge 2$, 
	for every $ t \in \{0, 1, \cdots, T\}$, $K_t: \Om \to 2^{\R^{d}}$ is a $\Fc^{0}_t$-measurable random set in the sense that
	$\{\om \in \Om : K_t(\om) \cap O \neq \emptyset \}  \in \Fc^{0}_t$ for every closed (open) set $O \subset \R^{d}$.
	Here, for each $\om \in \Om$, $K_t(\om)$ is a closed convex cone containing $\R^{d}_+$, called the solvency cone at time $t$. It represents the collection of positions, labelled in units of different $d$ financial assets, that can be turned into non-negative ones (component by component) by   performing immediately exchanges between the assets.  
	We denote by $K^*_t  \subset \R^{d}_+$ its (nonnegative) dual cone:
	\be \label{eq:def_K_star}
		K^*_t(\om) ~:=~ \big \{y \in \R^d  ~: \scap{ x}{y } \ge 0 ~\mbox{for all}~ x \in K_t(\om) \big \},
	\ee
	where $\scap{ x}{y } := \sum_{i=1}^{d} x^i y^i$ is the inner product on $\R^{d}$.
	For later use, let us also  introduce
	$$
		K^{*,0}_t(\om) ~:=~ \big\{ y = (y^1, \cdots, y^d) \in K^*_t(\om), ~y^d = 1\big \}.
	$$
	As in \cite{bouchard2016consistent}, we assume the following conditions throughout the paper:
		\begin{Assumption} 
	   $K^*_t \cap \partial \R^d_+ = \{0 \}$ and $\mathrm{int}K^*_t(\om) \neq \emptyset$ for every $\om \in \Omega$ and $t \le T$.
				
	\end{Assumption}
	It follows from the above assumption and   \cite[(viii)-Lemma A.1]{bouchard2016consistent} that there is a $\F^{0}$-adapted process $S$ satisfying
	\begin{align}\label{eq: S in int}
		S_t(\om) \in K^{*,0}_t(\om)\cap \mathrm{int}K^{*}_t(\om)\;\mbox{ for every $\omega \in \Om$, $t\le T$.}
	\end{align}
	We also assume that transaction costs are bounded and uniformly strictly positive. This is formulated in terms of $S$ above.   
	\begin{Assumption}\label{ass : assumptions}
		  There is some constant $c > 1$ such that
		$$
			c^{-1} S^i_t(\om) \le y^i \le c S^i_t(\om),
			~\mbox{for every}~i\le  d-1 
			~\mbox{and}~ y   \in K^{*,0}_t(\om).
		$$
	\end{Assumption}

	\begin{Example}\label{ex: example K}
		Let us consider a market with one risky asset with mid price $S^1_t > 0$ and one risk-free asset  $S^2_t \equiv 1$. Here $d=2$.
		Because of a proportional transaction cost parametrized by $c \ge 1$, 
		the bid price of the risky asset is given by $c^{-1} S^1_t$ and the ask price is $cS^1_t$.
		Then
		$$
			K_t(\om):=\{x\in \R^{2}: x^1 c^{-1} S^1_t(\om) \1_{\{x^1 \ge 0\}} + x^1 c S^1_t(\om) \1_{\{x^1 < 0\}} 
			+x^2 
			\ge 0\},
		$$
		$$
		 	K_t^*(\om) \!=\! \big\{(y^1, y^2) \in \R_+^2 : y^1 \!\in\! \big[y^2c^{-1} S^1_t(\om), y^2cS^1_t(\om) \big] \!\big\},
		$$
		and
		$$
		 	K_t^{*,0}(\om) \!=\! \big\{(y^1, 1) \in \R_+^2 : y^1 \!\in\! \big[c^{-1} S^1_t(\om), cS^1_t(\om) \big] \!\big\}.
		$$		
		Although there is, in the above example, a risk-free asset $S^2_t \equiv 1$ which serves as a num\'eraire,  this is not required in general. We refer to \cite{bouchard2016consistent} for an example with $d$ risky assets. 
		See also the monograph \cite{KabanovSafarian.09}.
	\end{Example}

Let us now turn to the definition of admissible trading strategies. 

	\begin{Definition} \label{def:admis_stra}
		We say that a $\F$-adapted process $\eta = (\eta_t)_{0 \le t \le T}$ is an admissible trading strategy if 
		$$ 
		\eta_t \in - K_t\;\;\;\Pc\mbox{-q.s.} \;\mbox{for all $t\le T$.}
		$$
		We denote by $\Ac$ the collection of all admissible strategies. 	
	\end{Definition}
The constraint $\eta_{t} \in -K_{t}$ means that $0-\eta_{t}\in K_{t}$, i.e., starting at $t$ with $0$, one can perform immediate transfers to reach the position $\eta_{t}$. 
Then, given $\eta\in \Ac$, the corresponding wealth process associated to a zero initial endowment at time $0$ is 
		$\big(\sum_{s=0}^{t} \eta_{s} \big)_{t\le T}$.

	\begin{Example} In the context of Example \ref{ex: example K}, $\eta\in \Ac$ if and only if 
		$$
		 \eta^1_{t} c^{-1} S^1_t  \1_{\{\eta_{t}^1 \le 0\}} + \eta_{t}^1 c S^1_t  \1_{\{\eta_{t}^1 > 0\}} 
			+\eta_{t}^2 
			\le 0\;\;\;\; \Pc\mbox{-q.s.} \;\mbox{for all $t\le T$.}
 		$$ 
	\end{Example}

\subsection{The randomization  approach}

	As explained in the introduction, we aim at considering a frictionless market set on an enlarged probability space, that is equivalent (in a certain sense) to our original market. This will be used later on to apply results  that are  already known  in the frictionless setting.  

	\vspace{1mm}

	Let us therefore first introduce an enlarged space.
	Let $c > 1$ be the constant in Assumption \ref{ass : assumptions}, we define $\Lambda_1 := [c^{-1}, c]^{d-1}$, $\Lambda_{t}:=(\Lambda_{1})^{t+1}$, and $\Lambda := \Lambda_T$,
	and then introduce the canonical process $\Theta_t(\theta) := \theta_t$, $\forall \theta = (\theta_t)_{0\le t \le T}  \in \Lambda$,
	as well as the $\sigma$-fields $\Fc_t^{\Lambda} := \sigma (\Theta_s,~ s \le t)$, $t\le T$.
	We next introduce an enlarged space $\Omb := \Om \x \Lambda$,
	an enlarged $\sigma$-field $\Fcb:= \Fc \otimes \Fc_T^{\Lambda}$, 
	together with two filtrations $\Fbb^0 = (\Fcb^0_t)_{0 \le t \le T}$ and $\Fbb = (\Fcb_t)_{0 \le t \le T}$ in which  $\Fcb^0_t := \Fc_t \otimes \{\emptyset, \Lambda\}$ and $\Fcb_t := \Fc_t \otimes \Fc_t^{\Lambda}$ for $t\le T$.

	Then, we define our randomized fictitious market by letting the fictitious stock price   $X = (X_t)_{0 \le t \le T}$ be defined by 
	\be \label{eq:def_X_t}
		X_t (\omb) := \Pi_{K^{*,0}_t(\om)} [S_t(\om) \theta_t], ~~\mbox{for all}~\omb = (\om, \theta) \in \Omb,\; t\le T, 
	\ee
	where $S_t(\om)\theta_t := (S_t^1(\om)\theta_t^1, \cdots, S_t^{d-1}(\om) \theta_t^{d-1}, S_t^d(\om))$,
	and $\Pi_{K^{*,0}_t(\om)}[y]$ stands for the projection of $y \in \R^d$ on the convex closed set $K^{*,0}_t(\om)$. Recall that   $S_t\in K^{*,0}_t$ for $t\le T$. 
	Finally, we introduce
	$$
		\Pcb ~:=~ \big\{  \Pb\in \Bf(\Omb, \Fcb) ~\mbox{such that}~\Pb|_{\Om} \in \Pc \big\}.
	$$

\begin{Remark}\label{rem: equiv pour tout et bar P} 
	We shall use several times the following important property related to the structure of $\Pcb$. Let $Y:(\om,\theta)\in \Omb\mapsto Y(\om,\theta)\in \R$ be a  random variable. Then, the following are equivalent:
\begin{itemize}
\item[\rm (i)] $Y\ge 0$  $\Pcb$-q.s
\item[\rm (ii)]  $Y(\cdot,\theta)\ge 0$ for all $\theta \in \Lambda$  $\Pc$-q.s.
\end{itemize}
The fact that {\rm (ii)} implies {\rm (i)} is clear. As for the reverse implication, we observe that if $\P\in \Pc$ is such that 
$\P[\inf_{\theta\in \Lambda} Y(\cdot,\theta)< -c]>0$ for some $c>0$, then one can find a Borel map $\om\in \Om\mapsto \theta(\om)\in \Lambda$ such that   $\P[  Y(\cdot,\theta(\cdot))< -c/2]>0$, see \cite[Lemma 7.27 and Proposition 7.49]{BertsekasShreve.78}. Since $\P\otimes \delta_{\theta(\cdot)}\in \Pcb$, this contradicts {\rm (i)}. \\
Moreover,  if $\theta\mapsto Y(\cdot,\theta)$ is upper semicontinuous, then  {\rm (i)}, {\rm (ii)} are equivalent to 
\begin{itemize}
\item[\rm (iii)]  $Y(\cdot,\theta)\ge 0$  $\Pc$-q.s.~for all $\theta \in \Lambda$.
\end{itemize}
\end{Remark}

	If one keeps the filtration $\F$ (or equivalently $\Fbb^0$) to define admissible (or self-financing) trading strategies on this fictitious financial market, then they coincide  with the admissible strategies in the sense of Definition \ref{def:admis_stra}  above. More precisely, we have the following. 

	\begin{Theorem} \label{thm:main}
		\rmi Fix $t \le  T$ and $\zeta_t\in L^{0}(\Fc_{t})$. Then
			$$
				\zeta_t \in -K_t~ \Pc \mbox{-q.s.}
				~~~\mbox{if and only if}~~~
				\scap{ \zeta_t}{ X_t }  \le 0~\Pcb \mbox{-q.s.}
			$$
		\noindent \rmii Consequently, an $\F$-adapted process $\eta$  is an admissible strategy in the sense of Definition \ref{def:admis_stra} 
		if and only if  $\scap{ \eta_t}{ X_t } \le 0$ $\Pcb$-q.s.~for all $t\le  T$.
	\end{Theorem}
	\proof The assertion (ii) is an immediate consequence of (i). To show that (i) holds, let us first note that,  
	by definition,  
	\be \label{eq:caract_K_t}
		x \in -K_t(\om)
		~\Longleftrightarrow~
		\scap{ x}{ y } \le 0~\forall y \in K_t^{*,0}(\om).
	\ee
	\noindent (a) First assume that $\zeta_t \in -K_t, \Pc$-q.s.~and fix  $\Pb \in \Pcb$ together with $\P := \Pb|_{\Om}$.  Then,  $\zeta_{t}(\om) \in -K_t(\om)$ for $\P$-a.e. $\om\in \Omega$.
	Recalling    \eqref{eq:caract_K_t}, this implies that 
	$\scap{ \zeta_{t}(\om)}{X_t(\om, \theta) } \le 0$, for every $\theta \in \Lambda$ and $\P$-a.e. $\om \in \Om$.
	Hence, $\scap{ \zeta_{t}}{ X_t } \le 0$ $\Pb$-a.s. by Remark  \ref{rem: equiv pour tout et bar P}.  By arbitrariness of $\Pb$, the later holds $\Pcb$-q.s.
	
	\vspace{1mm}
	
	\noindent (b) We now assume that $\scap{ \zeta_{t}}{X_t } \le 0~\Pcb$-q.s.  	
	Let $\P \in \Pc$, then for every $\theta \in \Lambda$, one has $\P \otimes \delta_{\theta} \in \Pcb$.
	Let $\Lambda_{\circ} \subset \Lambda$ be a countable dense subset,
	it follows that, for $\P$-a.e. $\om$ and every $y \in \{ S_{t}(\om) \theta_t ~:\theta\in \Lambda_{\circ}\}$, 
	one has $\scap{ \zeta_t(\om)}{y } \le 0$.
	By continuity of the inner product and Assumption \ref{ass : assumptions}, we then have 	  $\scap{ \zeta_{t}}{y } \le 0$ for all $y\in    K^{*,0}_{t}$ $\P$-a.s. The measure $\P\in \Pc$ being arbitrary chosen, we then deduce from  \eqref{eq:caract_K_t} that $\zeta_t \in -K_t, \Pc$-q.s. \endproof

%

\subsection{Equivalence of the no-arbitrage conditions under the framework of Bouchard \& Nutz \cite{bouchard2016consistent}}
\label{subsec:frameworkBN14}

	For the newly introduced fictitious market under model uncertainty, 
	a first issue is to formulate a no-arbitrage condition.
	Let us now specialize  to the framework of Bouchard \& Nutz \cite{bouchard2016consistent}.
	In this probabilistic framework, we show that the quasi-sure no-arbitrage condition of second kind used in \cite{bouchard2016consistent} is equivalent to the 
	quasi-sure no-arbitrage condition of \cite{BouchardNutz.13} on the frictionless market defined on $\Omb$ with stock price process $X$ and $\Fbb$-predictable strategies.
		
\paragraph{No-arbitrage condition under Bouchard \& Nutz's \cite{bouchard2016consistent} framework}	
	We first recall the framework of Bouchard \& Nutz \cite{bouchard2016consistent}. 
	Let  $\Om_0 =\{\om_0\}$ be a singleton and $\Om_1$ be a Polish space.
	For each $t \in \{1, \cdots, T\}$, we denote by $\Om_t := \Om_0 \x \Om_1^t$,
	where $\Om_1^t$ denotes the $t$-fold Cartesian product of $\Om_1$, we set 
	$\Fc^0_t := \Bc(\Om_t)$ and let $\Fc_t$ be its universal completion.	
	In particular, the $\sigma$-field $\Fc^0_0$ and $\Fc_0$ are trivial.
	From now on, 
	$$
		\Om := \Om_T, \;\;
		\Fc := \Fc_T,\;\;
		\F^0 := (\Fc^0_t)_{0 \le t \le T},\;\;
		\F := (\Fc_t)_{0 \le t \le T}.
	$$
	Given $t \in \{0, \cdots, T-1\}$ and $\om \in \Om_t$,   we are given a  non-empty {convex} set $\Pc_t(\om) \subseteq \Bf(\Om_1)$, which represents the set of possible models for the $(t+1)$-th period, given state $\om$ at time $t$.
	We assume that for each $t$ 
	\be \label{eq:AnalyticGraph}
		\mbox{graph}(\Pc_t)
		:= 
		\{ (\om, \P): \om \in \Om_t, 
			\P \in \Pc_t(\om) 
		\} 
		~\subseteq~ \Om_t \times \Pc(\Om_1)
		~~ \text{is analytic.}
	\ee
	Given   probability  kernels $\P_t:\Omega_{t}\mapsto \Bf(\Omega_{1})$,  for each $t \le T-1$, 
	we define a probability measure $\P$ on $\Om$ by Fubini's theorem:
	$$
 		\P(A) 
		:=
		\int_{\Om_1} \cdots \int_{\Om_1} \mathbf{1}_A (\om_1, \om_2 \cdots, \om_T) 
		\P_{T-1}(\om_1, \cdots, \om_{T-1}; d \om_T)  \cdots \P_0(d \om_1)
	$$
	We can then introduce the set $\Pc \subseteq \Bf(\Om)$ of possible models for the multi-period market up to time $T$:
	\be \label{eq:def_Pc}
		\Pc := 
		\big\{
			\P_0 \otimes \P_1 \otimes \cdots \otimes \P_{T-1} ~: \P_t(\cdot) \in \Pc_t(\cdot) \mbox{ for } t\le T-1 
		\big\}.
	\ee
	Notice that the condition \eqref{eq:AnalyticGraph} ensures that $\Pc_t$ admits  a universally measurable selector: $\P_t : \Om_t \rightarrow \Pc(\Om_1)$ such that $\P_t(\om) \in \Pc_t(\om)$ for all $\om \in \Om_t$.
	Then the set $\Pc$ defined in \eqref{eq:def_Pc} is nonempty.
	
	 \vspace{2mm}
 
	 Let us now recall the no-arbitrage condition used in \cite{bouchard2016consistent}.
	\begin{Definition}
		We say that $\NA2(\Pc)$ holds if for all $t \le  T-1 $ and all $\zeta \in L^0(\Fc_t)$, 
		\b*
			\zeta \in K_{t+1} ~~ \Pc \mbox{-q.s.} ~~~ \mbox{implies} ~~~ \zeta \in K_{t} ~~ \Pc \mbox{-q.s.}
		\e*
	\end{Definition}
	The following robust version of the fundamental theorem has been proved in \cite{bouchard2016consistent}. 
	\begin{Theorem} \label{thm:FundThm}
		The condition $\NA2(\Pc)$ is equivalent to :
		For all $t \le  T-1 $, $\P \in \Pc$ and $Y \in L_{\Pc}^0(\Fc_t,{\rm int}K_t^{*})$, there exists $\Q \in \Bf(\Om)$ and a $\F^0$-adapted process $(Z_s)_{s=t, \ldots, T}$ such that 
		$\P \ll \Q$ and $\P = \Q$ on $\Fc_t$, and

		\noindent \rmi $\Q \lll \Pc$

		\noindent \rmii  $Y=Z_t ~~ \Q\mbox{-a.s.}$

		\noindent \rmiii $Z_s \in {\rm int}K_s^* ~\Q\mbox{-a.s.~ for~} s=t, \ldots, T$

		\noindent \rmiv $(Z_s)_{s=t, \ldots, T}$ is a $\Q$-martingale, i.e.~$\E^{\Q}[Z_{s'}|\Fc_{s}]=Z_{s}$ for $t\le s\le s'$. 
	\end{Theorem}
	A  couple $(\Q, Z)$ satisfying the  conditions $\mathrm{(i)-(iv)}$ above for $t=0$ is called a strictly consistent price system (SCPS).  	
	For later use, let  $\CPSt$ denote the collection of all SCPS, and set   
		\be \label{eq:CPS0}
		\CPSt_0
		~:=
		\big\{ 
		(\Q, Z)\in \CPSt
		~\mbox{such that}~
		Z^d \equiv 1
		\big\}.
	\ee
	
	 For later use, we also recall the notion of $\NA2(t,\om)$ for each $t \le T$ and $\om \in \Om_t$:
	we say $\NA2(t, \om)$ holds true if 
	\be \label{eq:NA2t}
		\zeta \in K_{t+1}(\om, \cdot)~ \Pc_t(\om) \mbox{-q.s.}
		~~~\mbox{implies}~~
		\zeta \in K_t(\om), ~~\mbox{for all}~\zeta \in \R^d.
	\ee
	The following result is proved in  \cite[Lemma 3.6]{bouchard2016consistent}.
	\begin{Lemma} \label{lemm:NA2_polar}
		The set $N_t := \{\om : \NA2(t,\om) ~\mbox{fails}\}$ is universally measurable.
		Moreover,  $N_t$ is a $\Pc$-polar set if $\NA2(\Pc)$ holds.
	\end{Lemma}

\paragraph{No-arbitrage condition on the enlarged space}
	We next consider the enlarged space $(\Omb, \Fcb)$ and define a subset of probability measures 
	$\Pcb_{\mathrm{int}} \subset \Pcb$,
	in order to introduce the quasi-sure no-arbitrage condition  of \cite{BouchardNutz.13} w.r.t.~the price process $X$ and the  set of strategies
	$$\Hcb := \{ \mbox{All}~ \Fbb \mbox{-predictable processes} \}.$$
	Given $t \le T$, we denote 
	$\Omb_0 := \Om_0 \x \Lambda_1$ and $\Omb_t := \Omb_0 \x (\Om_1 \times \Lambda_1)^t$, so that one has
	$\Omb := \Om \x \Lambda = \Omb_T$. We write $(\om,\theta)\in \Omb$ in the form $\om=(\om_{0},\ldots,\om_{T})$ and  $\theta=(\theta_{0},\ldots,\theta_{T})$. For $t\le T$, we use the notations $\omega^{t}:=(\om_{0},\ldots,\omega_{t})$, $\theta^{t}=(\theta_{0},\dots,\theta^{t})$ and $\omb^{t}=(\om^{t},\theta^{t})$. 
	We now introduce the subset $\Pcb_{\mathrm{int}} \subset \Pcb := \{\Pb \in \Bf(\Omb) ~: \Pb|_{\Om} \in \Pc \}$ defined as follows:
	\begin{itemize}	
		\item For $t = 0, 1, \cdots, T-1$ and $\omb = (\om, \theta) \in \Omb_t$,
		define $\Pcb_t(\omb) := \big\{ \Pb \in \Bf(\Om_1 \x \Lambda_1) ~: \Pb|_{\Om_1} \in \Pc_t(\om) \big\}$, and 
		\be \label{eq:def_Pcb_t}
			\Pcb_t^{\mathrm{int}}(\omb) := \big\{\Pb \in \Pcb_t(\omb) ~: (\delta_{\omb^{t}} \otimes \Pb) [X_{t+1} \in \mathrm{int} K^*_{t+1}] = 1 \big\},
		\ee
		where $\delta_{\omb^{t}} \otimes \Pb$ is a probability measure on $\Omb_{t+1} = \Omb_t \x (\Om_1 \x \Lambda_1)$ and 
		$X_{t+1}$ (defined in \eqref{eq:def_X_t}) is considered as a random variable defined on $\Omb_{t+1}$.

		\item Let $\Pcb^{\mathrm{int}}_{\emptyset}$ be the collection of all probability measures $\Pb$ on $\Omb_0$ such that $\Pb[X_0 \in \mathrm{int} K^*_0] = 1$, and define
		$$
			\Pcb_{\mathrm{int}}
			~:=~
			\big\{
				\Pb_{\emptyset} \otimes \Pb_0 \otimes \cdots \otimes \Pb_{T-1} ~: 
				\Pb_{\emptyset} \in \Pcb^{\mathrm{int}}_{\emptyset}
				~\mbox{and}~
				\Pb_t(\cdot) \in \Pcb^{\mathrm{int}}_t(\cdot) ~\mbox{for}~ t \le T-1
			\big\},
		$$
		where $\Pb_t(\cdot)$ is a universally measurable selector of $\Pcb^{\mathrm{int}}_t(\cdot)$, whose existence is ensured by
		Lemma \ref{lemma:panalytic} below.
	\end{itemize}
	
	 By a slight abuse of notations, we shall later write $\Pcb_t(\omb)$ for $\Pcb_t(\omb^t)$ when 
	$\omb \in \Omb$.
	The same convention will be used for $\Pcb^{\mathrm{int}}_t(\omb)$, etc.

	\begin{Remark}\label{rem: equiv pour tout et bar P bis} 
		Note that the equivalence observed in Remark \ref{rem: equiv pour tout et bar P} still holds for the above construction. In particular, let $Y:(\om,\theta)\in \Omb\mapsto Y(\om,\theta)$ be a random variable. Then, the following are equivalent:
\begin{itemize}
\item[\rm (i)] $Y\ge 0$  $\Pcb_{\mathrm{int}}$-q.s
\item[\rm (ii)]  $Y(\om,\theta)\ge 0$ for all $\theta \in \Lambda^{\mathrm{int}}(\om) := \{\theta' :S_{t}(\om)\theta'_{t}\in {\rm int}K^{*}_{t}(\om), ~\forall t\le T\}$,  for all $\om$ outside a $\Pc$-polar set.
\end{itemize}
	\end{Remark}
	
	\begin{Lemma} \label{lemma:panalytic}
		For every $0 \le t \le T-1$,the set
		\b*
			\mathrm{graph} \big(\Pcb_t^{\mathrm{int}} \big)
			:= 
			\big\{ (\omb, \Pb): \omb \in \Omb_t, 
				\Pb \in \Pcb^{\mathrm{int}}_t(\omb) 
			\big\} 
			~\text{is analytic.}
		\e*
	\end{Lemma}
	\proof We only consider the case $t \ge 1$,  the proof for the case $t=0$ is an obvious modification.
	
	\rmi Since $\mathrm{graph}(\Pc_t)$ is an analytic set, there is some Polish space $E$ and a Borel set
	$A \subset \Om_t \x \Bf(\Om_1) \x E$ such that $\mathrm{graph}(\Pc_t)$ is the projection set of $A$ on $\Om_t \x \Bf(\Om_1)$, i.e. $\mathrm{graph} \big(\Pc_t \big) = \Pi_{\Om_t \x \Bf(\Om_1)} \big[ A \big]$.
	Let us define
	$$
		\overline A 
		~:=~ 
		\big\{ \big(\omb, \Pb, e \big)  \in \Omb_t \x \Bf(\Om_1 \x \Lambda_1) \x E~: \big(\om, \Pb|_{\Om_1}, e\big) \in A
		\big\},
	$$
	which is a Borel set in $\Omb_t \x \Bf (\Om_1 \x \Lambda_1) \x E$ since $\omb = (\om, \theta) \mapsto \om$ and $\Pb \mapsto \Pb|_{\Om_1}$ are Borel.
	Then $\mathrm{graph}(\Pcb_t) = \Pi_{\Omb_t \x \Bf(\Om_1 \x \Lambda_1)} \big[~\overline A~ \big]$ is an analytic set.
	
	\rmii By Lemma A.1 of \cite{bouchard2016consistent}, we know that for every $t \ge 1$, 
	$\{(\om^{t+1}, \theta^{t+1}, x) \in \Omb_{t+1} \x \R^d ~: x \in \mathrm{int} K^*_{t+1}(\om^{t+1}) \}$ is Borel measurable.
	Then 
	$$ (\om^t, \theta^t, \Pb) \mapsto \big( \delta_{(\om^t, \theta^t)} \otimes \Pb \big) \big[X_{t+1} \in \mathrm{int} K^*_{t+1} \big]$$
	is also Borel measurable and hence
	$$
		B := \Big\{ \big(\omb^t, \Pb \big) \in \Omb_t \x \Bf(\Om_1 \x \Lambda_1) ~: 
		\big(\delta_{\omb^t} \otimes \Pb\big) \big[X_{t+1} \in \mathrm{int} K^*_{t+1} \big] = 1 \Big\}
	$$
	is a Borel set.
	Then it follows that $\mathrm{graph} \big(\Pcb_t^{\mathrm{int}} \big) = B \cap \mathrm{graph}(\Pcb_t) $ is an analytic set.
	\qed

	\begin{Definition}\label{def: NA}  We   say that $\NA_{t}(\Pcb_{\mathrm{int}})$ holds   if
		$$
			(\Hb \circ X)_T \ge 0,~ \Pcb_{\mathrm{int}}\mbox{-q.s.}
			~~~\Longrightarrow~~~
			(\Hb \circ X)_T= 0,~ \Pcb_{\mathrm{int}}\mbox{-q.s.},
		$$
		 for every  $\Hb\in  \Hcb$.   
	\end{Definition}
	Using \cite{BouchardNutz.13} and Lemma \ref{lemma:panalytic} above, it follows that the following fundamental theorem of asset pricing holds. 
\begin{Theorem} \label{thm:FundThmNT}
		The condition $\NA(\Pcb_{\mathrm{int}})$ is equivalent to :
		For all  $\Pb \in \Pcb_{\mathrm{int}}$, there exists $\Qb \in \Bf(\Omb)$ such that $\Pb \ll\Qb \lll \Pcb_{\mathrm{int}}$ and $X$ is a $(\Qb, \Fbb)$-martingale.
	\end{Theorem}
	Hereafter, we denote by $\Qcb_{0}$ the collection of measures $\Qb \in \Bf(\Omb)$ such that $\Qb \lll \Pcb_{\mathrm{int}}$ and $X$ is a $(\Fbb, \Qb)$-martingale.

	\vspace{2mm}

	The main result of this section says that the two no-arbitrage conditions defined above are equivalent. 
	\begin{Proposition} \label{prop:equiv_NA}
		The conditions $\NA2(\Pc)$ and  $\NA(\Pcb_{\mathrm{int}})$ are equivalent.
	\end{Proposition}
	\proof (i) Let us first suppose that  $\NA(\Pcb_{\mathrm{int}})$ holds.
	Assume that $\NA2(\Pc)$ does not hold. Then,  for some $t \le T-1$, there is $\zeta \in L^0(\Fc_t)$ such that
	$\zeta \in K_{t+1}$ $\Pc$-q.s.~and $\P[A]>0$ for $A:=\{\zeta \notin K_{t}\}$ and some $\P\in \Pc$.
	Using a standard measurable selection argument and \eqref{eq: S in int},
	there is a measurable map $f:\R^{2d}\mapsto \Lambda_1$ such that 
	$Y_{t}:=S_{t}f(\zeta,S_{t}) \in L_\Pc^0(\Fc_t, K_t^{*,0}\cap \mbox{int}K_t^{*})$ and  $A\subset \{\scap{ Y_{t}}{\zeta}<0\}$. 
	Set $\bar \zeta:=\zeta\1_{\{\scap{ \zeta}{ X_{t}}\le 0\}}$,
	so that $\scap{\bar \zeta }{\big( X_{t+1}-X_{t} \big)}\ge 0$ $\Pcb_{\mathrm{int}}$-q.s.
	However, using Assumption \ref{ass : assumptions}, \eqref{eq: S in int} and the fact that $S$ does not depend on $\theta \in \Lambda$ while $A\in \Fc$,
	it follows that there is some $\Pb  \in \Pcb_{\mathrm{int}}$ such that 
	$\Pb[(X_{t},X_{t+1})=(Y_{t},S_{t+1})] = 1$
	and $\Pb[\scap{\bar \zeta }{(X_{t+1}-X_{t})}>0] \ge \P[A] >0$.
	This contradicts $\NA(\Pcb_{\mathrm{int}})$.  

	\vspace{1mm}

	\noindent (ii) Conversely, assume that $\NA2(\Pc)$ holds, we aim at proving that $\NA(\Pcb_{\mathrm{int}})$ holds.
	In view of Theorem \ref{thm:FundThmNT}, it is enough to prove that, for every $\Pb \in \Pcb_{\mathrm{int}}$, 
	there is a $\Qb \lll \Pcb_{\mathrm{int}}$ such that $\Pb \ll \Qb$ and $X$ is a $(\Fbb, \Qb)$-martingale.

	Fix $\Pb \in \Pcb_{\mathrm{int}}$, then, by the definition of $\Pcb_{\mathrm{int}}$, one has the representation:
	$$
		\Pb ~:=~ \Pb_{\emptyset} \otimes (\P_0 \otimes q_0) \otimes \cdots \otimes (\P_{T-1} \otimes q_{T-1}),
	$$
	where $\P_t: \Omb_t \to \Pc(\Om_1)$ and $q_t \big(d \theta_{t+1} | \omb^t, \om_{t+1} \big) : \Omb_t \x \Om_1 \to \Pc(\Lambda_1)$ are all Borel kernels such that
	$\P_t = (\P_{t}(\cdot |\omb^t))_{\omb^t \in \Omb_t}$ satisfies $\P_{t}(\cdot |\omb^t) \in \Pc(t, \om^t)$
	and  the support of $q_t(\cdot | \omb^t, \om_{t+1})$ is contained in 
	$\Lambda^{\mathrm{int}}_{t+1}(\om^{t+1})$ $:= $ $\{\theta_{t+1} \in \Lambda_1 ~: S_{t+1}(\om^{t+1}) \theta_{t+1} \in \mathrm{int}K^*_{t+1}(\om^{t+1}) \}$.
	
	We next construct another kernel $q'_t$ in order to define a martingale measure $\Qb$ dominating $\Pb$. For every $t \le T-1$, we consider the Borel kernel  $(\P_{t}(\cdot |\omb^t))_{\omb^t \in \Omb_t}$.
	It follows\footnote{It suffices to replace their $\Omega_{t}$ by our $\Omb_{t}$, and take their $Y$ equal to $X_{t}$.} from  \cite[Lemmas 3.8 and 3.9]{bouchard2016consistent}  that there is a universally measurable map $(\omb^t, \om_{t+1}) \mapsto Z_{t+1}(\omb^t, \om_{t+1})$ and two families $(\Q_{t}(\cdot |\omb^t))_{\omb^t \in \Omb_t}$ and  $(\P'_{t}(\cdot |\omb^t))_{\omb^t \in \Omb_t}$ such that
	$\P'_{t}(\cdot |\omb^t) \in \Pc(t, \om^t)$, $Z_{t+1} \in \mathrm{int} K^*_{t+1} \cap K^{*,0}_{t+1}$ and
	\be \label{eq:claim_NA}
		\P_{t}(\cdot |\omb^t) \ll \Q_{t}(\cdot |\omb^t) \ll \P'_{t}(\cdot |\omb^t),
		~~\mbox{and}~~
		\E^{\Q_{t}(\cdot |\omb^t)}[ Z_{t+1}] = S_t(\om^t) \theta_t = X_t(\omb^t).
	\ee
	We next consider the family $\big(q_t(\cdot| \omb^t, \om_{t+1}) \big)_{\omb^t \in \Omb_t, \om_{t+1} \in \Om_1}$.
	By  \cite[Lemma 4.7]{BouchardNutz.13}, the set of all $(\omb^t, \om_{t+1}, \alpha, q'(d \theta_{t+1}))$ satisfying
	$$
		q'[\Lambda^{\mathrm{int}}_{t+1}(\om^{t+1}) ] = 1,
		~~
		q' \gg q_t( \cdot | \omb^t, \om_{t+1})  
		~~\mbox{and}~
		\E^{q'} \big[X_{t+1}(\om^{t+1}, \cdot) \big] = \alpha,
	$$
	is a Borel set. It is not difficult to see that this set is non-empty.
	Then, by a standard measurable selection argument, there is a universally measurable family 
	$q'_t( \cdot| \omb^t, \om_{t+1})$  with support in  $\Lambda^{\mathrm{int}}_{t+1}(\om^{t+1})$ such that
	$$
		q_t( \cdot| \omb^t, \om_{t+1}) \ll q'_t( \cdot| \omb^t, \om_{t+1})
		~~~\mbox{and}~~
		\E^{q'_t( \cdot| \omb^t, \om_{t+1})}[X_{t+1}(\om^{t+1}, \cdot)] = Z_{t+1}(\omb^t, \om_{t+1}).
	$$
 
	Let us finally define
	$$
		\Qb := \Pb_{\emptyset} \otimes \big( \Q_0 \otimes q'_0 \big) \otimes \cdots \otimes \big( \Q_{T-1} \otimes q'_{T-1} \big)
	$$
	and
	$$
		\Pb' := \Pb_{\emptyset} \otimes \big( \P'_0 \otimes q'_0 \big) \otimes \cdots \otimes \big( \P'_{T-1} \otimes q'_{T-1} \big).
	$$
	Then it is easy to check that 
	$$
		 \Pb \ll \Qb \ll \Pb',
		~~
		\Pb' \in \Pcb,
		~~
		\Qb \in \Qcb_0,
	$$
	and we hence conclude the proof.
	\qed

	\begin{Remark}\label{rem: Na loc condi}	
		Let us define $\Lambda^{\rm int}_{0}(\om_{0}) := \{\theta_{0} \in \Lambda_{1}: S_{0}(\om_0) \theta_{0}\in {\rm int} K^{*}_{0} \}$, and for each $\theta_{0} \in  \Lambda_{0}^{\rm int}(\om_{0})$, 
		$$
			\Pcb^{\mathrm{int},\delta}_{0}(\theta_{0}) 
			\!:=\! 
			\big\{ \Pb \in \Pcb_{\mathrm{int}} ~: \Pb[ \Theta_0 = \theta_{0} ] = 1 \big\},
			~~~
			\Qcb^{\delta}_0(\theta_{0}) 
			\!:=\! 
			\big\{ \Qb \in \Qcb_0 ~: \Qb[\Theta_0 = \theta_{0} ] = 1 \big \}.
		$$
		 Define $\NA(\Pcb^{\mathrm{int},\delta}_{0}(\theta_{0}))$ as $\NA(\Pcb_{\mathrm{int}})$ with $\Pcb^{\mathrm{int},\delta}_{0}(\theta_{0})$ in place of  $\Pcb_{\mathrm{int}}$.
		Then,  
		$\NA(\Pcb_{\mathrm{int}})$ implies that $\NA(\Pcb^{\mathrm{int},\delta}_{0}(\theta_{0}))$ holds for  every
		$\theta_{0}\in \Lambda^{\rm int}_{0}(\om_{0})$. Indeed,  assume that $\NA(\Pcb_{\mathrm{int}})$ holds. Then,   Theorem \ref{thm:FundThmNT}  applied to $\Pcb_{\mathrm{int}}$ implies that, for any  $\Pb \in \Pcb_{\mathrm{int}}$ such that $\Pb[\Theta_0 = \theta_{0}] = 1$, one can find $\Pb   \ll \Qb(\theta_{0})$ such that $X$ is a $\Qb(\theta_{0})$-martingale. 
	\end{Remark}

\section{A robust pricing-hedging duality result}	
\label{sec:pricing_hedging}

	We now concentrate on the super-replication problem
	under the framework of \cite{bouchard2016consistent}. 
	Given $e \in \N\cup\{0\}$, we are given a random vector $\xi: \Om \to \R^{d}$ 
	as well as $\zeta_i: \Om \to \R^{d}$ for $i=1, \cdots, e$ such that $\zeta_i \not\equiv 0$.
	The random vectors $\xi$ and $\zeta_i$ represents the final payoffs, in number of units of each risky assets, of respectively an exotic option and vanilla options.
	We assume that the bid and ask prices of each vanilla option $\zeta_i$ are respectively $-c_i$ and $c_i$ for some constant $c_i  \ge  0$, this symmetry is without loss of generality.
	Then,  the minimal super-hedging cost of the exotic option $\xi$ using vanilla options $\zeta_i$ together with dynamic trading strategy is given by\footnote{Here we use the convention $\sum_{i=1}^{0} = 0$}:
	\begin{equation} \label{eq:super_hedging_e}
		\pi_e(\xi):=
		\inf \Big\{ y+ \sum_{i=1}^e c_i |\ell_i|~: 
			 y{\rm 1}_{d} + \sum_{i=1}^e \ell_i \zeta_i + \sum_{t =0}^T \eta_t - \xi \in K_T,~\Pc \mbox{-q.s.}, (\eta, \ell)\in \Ac\x \R^{e}
		\Big\},
	\end{equation}
	where ${\rm 1}_{d}$ is the vector will all components equal to $0$ but the last one that is equal to $1$. 
	Let us introduce the  subset of the set of  SCPS (recall \eqref{eq:CPS0}) that are compatible with the bid-ask speads of the vanilla options used for static hedging: 
	$$
		\CPSt_e := 
		\big\{ (\Q, Z) \in \CPSt_0
			~:
			\E^{\Q} \big[ \scap{ \zeta_i }{ Z_T } \big] \in [-c_i, c_i],~ i=1,\cdots, e
		\big\}.
	$$
	Then, we have the following super-hedging duality. 
	\begin{Theorem} \label{thm:main_surrep}
		Let $\xi$ and $(\zeta_i)_{i\le e}$ be   Borel measurable, and assume that  $\NA2(\Pc)$ holds true.
		Assume either that $e=0$, 
		or that $e \ge 1$  and 
		\be \label{eq: NA avec option et CT} 
			\sum_{i=1}^{e}   \big(\ell_i  \zeta_i - |\ell_{i}| c_{i} {\rm 1}_{d} \big)+\sum_{t=0}^{T} \eta_{t}\in K_{T}~\Pc\mbox{-q.s.}
			~~\Longrightarrow~~
			\ell =0
		\ee		
		for all $\ell \in \R^{e}$ and $\eta\in \Ac$.
		Then $\CPSt_e$ is nonempty and  
		\be \label{eq:pricinghedingduality}
			\pi_e(\xi) ~= \sup_{(\Q, Z) \in \CPSt_e} \E^{\Q} \big[ \xi \cdot Z_T \big].
		\ee
		Moreover, there exists $(\hat \eta, \hat \ell) \in \Ac \x \R^e$ such that 
		$$
			\pi_e(\xi){\rm 1}_{d} 
			+ \sum_{i=1}^e  \Big(  \hat \ell_i \zeta_i -  |\hat \ell_i| c_i {\rm 1}_{d} \Big) 
			+ \sum_{t =0}^T \hat \eta_t - \xi \in K_T,
			~\Pc \mbox{-q.s.}
		$$
	\end{Theorem}

	The proof is provided in the subsequent sections. We start with the case $e=0$.

\subsection{Proof of Theorem \ref{thm:main_surrep}: case $e=0$}
\label{subsec:proof_e0}

\subsubsection{Reformulation of the super-hedging problem}
\label{subsec: equiv problems}
	
	Before providing the proof of Theorem \ref{thm:main_surrep}, we first reformulate 
	the optimization problem  \eqref{eq:pricinghedingduality} and 
	the super-hedging problem \eqref{eq:super_hedging_e}, under proportional transaction cost, in terms of the   fictitious market defined on our enlarged space.

	We start with the pricing problem. 
	Let  us define
	$$
		\Qcb_0
		:=
		\big\{
			\Qb \lll \Pcb_{\mathrm{int}} ~: X ~\mbox{is a}~(\Fbb, \Qb) \mbox{-martingale s.t.}~X_t \in \mathrm{int} K^*_t~\forall t\le T
			~\Qb\mbox{-a.s.}
		\big\}
	$$
	and
	$$
		\Qcb^{\mathrm{loc}}_0
		:=
		\big\{
			\Qb \lll \Pcb_{\mathrm{int}} ~: X ~\mbox{is a}~(\Fbb, \Qb) \mbox{-local martingale s.t.}~X_t \in \mathrm{int} K^*_t~\forall t\le T
			~\Qb\mbox{-a.s.}
		\big\}.
	$$	
	
	\begin{Proposition} \label{prop:equiv_pricing}
		For any universally measurable vector $\xi: \Om \to \R^{d}$, one has
		$$
			\sup_{(\Q,Z) \in \CPSt_0} \E^{\Q} \big[\xi \cdot Z_T \big]
			~=~
			\sup_{\Qb \in \Qcb_0} \E^{\Qb} \big[ \xi \cdot X_T \big].
		$$
	\end{Proposition}
	\proof 
	\rmi 	
	First, let $(\Q, Z) \in \CPSt_0$ be a consistent price system, then there is a Borel measurable map $\rho_t: \Om_t \to \Lambda$ such that $X_t(\om, \rho_t(\om)) = S_t(\om) \rho_t(\om) = Z_t(\om)$ for all $\om \in \Om$.
	Let us define $\Qb$ by $\Qb[ \overline A] := \Q[ \{ (\om_{t}, \rho_t(\om))_{t\le T} \in \overline A \}]$ for every $\overline A \in \Fcb_T$.
	Then,  $ \Qb \in \Qcb_0$ and it satisfies
	$\E^{\Qb} \big[ \xi \cdot X_T \big] = \E^{\Q} \big[ \xi \cdot Z_T \big]$.

	Conversely, given $\Qb \in \Qcb_0$, let us  define $\Q := \Qb|_{\Om}$ and $Z_t := \E^{\Qb}[ X_t | \Fc_t]$ for $t\le T$.
	Since $\Qb \ll \Pb$ for some $\Pb \in \Pcb$, then $\Q \ll \P := \Pb |_{\Om} \in \Pc$.
	Moreover, the fact that $X$ is a $(\Fbb, \Qb)$-martingale implies that $Z$ is $(\F, \Q)$-martingale.
	Then,  $(\Q, Z)$ is a strictly consistent price system, and 
	$\E^{\Q} \big[ \xi \cdot Z_T \big] = \E^{\Qb} \big[ \xi \cdot X_T \big] $.
	\qed

	\vspace{2mm}

	We next reformulate the super-hedging problem \eqref{eq:super_hedging_e} on the enlarged space.
	Let us define
	$$
		g(\omb) := \scap{ \xi(\om) }{ X_T(\omb) },
		~\mbox{for all}~
		\omb = (\om, \theta) \in \Omb,
	$$
	as the contingent claim.  
	Denote by $\Hc$  the collection of all $\Fbb^0$-predictable processes
	and let 
	$(H \circ X)_t := \sum_{s = 1}^t \scap{ H_{s}}{( X_s - X_{s-1})} $, $t\le T$, 
	be the wealth process associated to $H \in \Hc$.

	\begin{Proposition} \label{prop:reformulation} 
		One has 
		\b*
			\pi_0(\xi) 
			&=&
			\inf \Big\{
				y  \in \R~: y + (H \circ X)_T \ge g ~\Pcb \mbox{-q.s.}, \mbox{ for some }H \in \Hc
			\Big\}\\
			&=&
			\inf \Big\{
				y  \in \R~: y + (H \circ X)_T \ge g ~\Pcb_{\mathrm{int}} \mbox{-q.s.}, \mbox{ for some }H \in \Hc
			\Big\}.
		\e*
	\end{Proposition}
	\proof For ease of notations, we write $\Delta X_{t}:=X_{t}-X_{t-1}$. 
	
	\noindent \rmi  Let $(y,\eta) \in \R\times\Ac$ be such that  and $y{\rm 1}_{d} + \sum_{t=0}^T \eta_t - \xi \in K_T$  $\Pc$-q.s.
	Define the $\Fbb^0$-predictable process $H$ by $H_t := \sum_{s = 1}^t \Delta H_s$ with
	$\Delta H_t := \eta_{t-1}$, for $t = 1, \cdots, T$.
	By exactly the same arguments as in part $\mathrm{(i)}$ of the proof of Theorem \ref{thm:main}, this is equivalent to 
	\begin{align} 
		0&\le \scap{ \Big( y{\rm 1}_{d} + \sum_{t=0}^{T} \eta_t - \xi \Big)}{ X_T } \nonumber
		\\
		&= y +  \sum_{t=0}^{T} \scap{ \eta_t}{( X_T-X_{t} )}+  \sum_{t=0}^{T} \scap{ \eta_t}{ X_t } - g \nonumber\\
		&= y + \sum_{t=1}^T \scap{ H_t}{\Delta X_t} + \sum_{t=0}^{T} \scap{ \eta_t}{ X_t } - g ~~~\Pcb \mbox{-q.s.},\label{eq:sup_rep_equiv}
	\end{align}
	where the last equivalence follows by   direct computation using that  $X^d_t \equiv 1$. Since $ \scap{ \eta_t}{  X_t } \le 0$ $\Pcb$-q.s., 
	 by Theorem \ref{thm:main}, we deduce that 
	$y + (H \circ X)_T \ge g ~\Pcb$-q.s. This shows that 
	$$
		\pi_0(\xi) 
		~\ge~
		\inf \big\{ y\in \R ~: y + (H \circ X)_T \ge g~\Pcb\mbox{-q.s.~for some }  H \in \Hc \big\}.
	$$
	
	\noindent \rmii We next prove the converse inequality. 
	Let $(y, H) \in \R \x \Hc$ be such that $y + (H \circ X)_T \ge g~\Pcb\mbox{-q.s.}$ We use the convention $H_{0}=0$.  Set $\eta^i_t := \Delta H^i_{t+1}$ for all $i=1, \cdots, d - 1$ and $t\le  T-1$, and $\eta_T := 0$. 
	We next define $\eta^d_t$ for $t= 0, \cdots, T-1$ by
	\be \label{eq:def_eta_d}
		\eta^d_t(\om) := \min_{\theta \in \Lambda} m^d_t(\om, \theta)
		~~\mbox{with}~ m^d_t(\omb) := - \sum_{i=1}^{d-1} \eta^i_t(\om) X^i_t(\omb),
	\ee
	for all $\omb = (\om, \theta) \in \Omb$.
	Notice that $m^d_t(\om, \theta)$ is bounded continuous in $\theta$,
	then $\eta^d_t(\om) = \min_{\theta \in \Lambda_\circ} m^d_t(\om, \theta)$ for some countable dense subset $\Lambda_\circ$ of $\Lambda$, and hence $\eta^d_t \in \Fc_t$.
	Using its construction, one has $\eta \in \Ac$. 
	Moreover, it follows from the choice of $(y,H)$ and  the fact that $\P \times \delta_{\theta} \in \Pcb$ for all $\P \in \Pc$ that    
	\begin{align}\label{eq: inf theta pour ega pi}
	0\le & \inf_{\theta \in \Lambda_\circ }  \big( y + (H \circ X)_T - g \big) (\cdot, \theta) 
	~=~ \inf_{\theta \in \Lambda }  \big( y + (H \circ X)_T - g \big) (\cdot, \theta) \\
	=& \inf_{\theta \in \Lambda }   \Big( y +  \sum_{t=0}^{T} \scap{ \eta_t}{( X_T-X_{t}) } - \scap{ \xi}{ X_{T}}  \Big) (\cdot, \theta)\nonumber\\
	= &  \inf_{\theta \in \Lambda }  \Big(   \scap{ \big( y{\rm 1}_{d}+\sum_{t=0}^{T} \eta_t-\xi \big)}{ X_T  }- \sum_{t=0}^{T-1} \scap{ \eta_t}{  X_{t} }    \Big) (\cdot, \theta)  ~\Pc \mbox{-q.s.},\nonumber
	\end{align}
	recall that $\eta_{T}=0$ by its construction above \eqref{eq:def_eta_d}. 
	We now use the definition of $\eta^{d}$  and the fact that each $X_{t}$ depends on $\theta$ only through $\theta_{t}$ to obtain 
	\begin{align*}
	0\le &     \inf_{\theta \in \Lambda }  \big\{\scap{ (y{\rm 1}_{d}+\sum_{t=0}^{T} \eta_t-\xi)}{ X_T  }\big\} (\cdot, \theta)- \sum_{t=0}^{T-1} \sup_{\theta \in \Lambda} \big\{\scap{ \eta_t}{ X_{t} } \big\} (\cdot, \theta) \\
	=& \inf_{\theta \in \Lambda }  \big\{\scap{ (y{\rm 1}_{d}+\sum_{t=0}^{T} \eta_t-\xi)}{ X_T  } \big\}(\cdot, \theta) ~\Pc \mbox{-q.s.}
	\end{align*}
	The latter is equivalent to  $y{\rm 1}_{d}+ \sum_{t=0}^T \eta_t - \xi \in K_T$, $\Pc$-q.s. This shows that 
	$$
		\pi_0(\xi) 
		~\le~
		\inf \big\{ y\in \R ~: y + (H \circ X)_T \ge g ~\Pcb\mbox{-q.s.},  \mbox{ for some }H \in \Hc \big\}.
	$$	
	\noindent \rmiii Let us now prove that 
	$$
		\pi_0(\xi)=\inf \big\{ y\in \R~: y + (H \circ X)_T \ge g ~\Pcb_{\mathrm{int}}\mbox{-q.s.},  \mbox{ for some }H \in \Hc \big\}.
	$$

	Since $\Pcb_{\mathrm{int}} \subset \Pcb$, one  inequality follows from $\mathrm{(i)-(ii)}$ above. 
	As for the converse one, 
	let $(y, H) \in \R \x \Hc$ be such that $y + (H \circ X)_T \ge g ~\Pcb_{\mathrm{int}}$-q.s.
	and   define $\eta$ as in \eqref{eq:def_eta_d}.
	Observe that 
	the right-hand side term of \eqref{eq: inf theta pour ega pi} is equal to 
	$$
	 \inf_{\theta \in \Lambda^{{\rm int}}(\cdot)}  \big( y + (H \circ X)_T - g \big) (\cdot, \theta)
	$$
	$\Pc$-q.s., in which, for $\omega\in \Omega$, $\Lambda^{{\rm int}}(\omega)$ is
	defined as the collection of $\theta \in \Lambda$ such that $S_{t}(\omega)\theta_{t}\in {\rm int}K^{*}_{t}(\omega)$ for all $t\le T$. 
	
 Next, to each $\theta \in \Lambda$,
	we associate the probability kernels 
	\be \label{eq:kernel_q_theta}
		q^{\theta}_{s}: \om \in \Omega \mapsto  q^{\theta}_{s}(\cdot|\om)
		:=
		\delta_{\theta}\1_{A^{\theta}_{s}(\om)} 
		+\delta_{{\bf 1}}\1_{(A^{\theta}_{s}(\om))^{c}}\in \Bf(\Lambda),~ s\le T,
	\ee
	where  ${\bf 1}$ is the vector of $\R^{d}$ with all entries equal to $1$, 
	$A^{\theta}_{s}(\om):=\emptyset$ for $s\ne t$ and $A^{\theta}_{t}(\om):=\{ S_{t}(\om) \theta_{t}\in {\rm int}K^{*}_{t}(\om)   \}$.
	It follows that $\P\otimes (q^{\theta}_{0}\otimes q^{\theta}_{1}\otimes\cdots\otimes q^{\theta}_{T}) \in \Pcb_{\mathrm{int}}$ for every $\P \in \Pc$.
	Then it suffices to argue as in \rmii above to obtain that
	\begin{align*}
	0&\le \inf_{\theta \in \Lambda^{{\rm int}}(\cdot)}  \big( y + (H \circ X)_T - g \big) (\cdot, \theta)=\inf_{\theta \in \Lambda}  \big( y + (H \circ X)_T - g \big) (\cdot, \theta) ~\Pc \mbox{-q.s.}
	\end{align*}
	which implies that
	\begin{align*}
	0\le & \inf_{\theta \in \Lambda}  \big\{\scap{ (y{\rm 1}_{d}+\sum_{t=0}^{T} \eta_t-\xi)}{ X_T  } \big\}(\cdot, \theta) ~\Pc \mbox{-q.s.,}
	\end{align*}
	and we hence conclude as in step $\mathrm{(ii)}$.
	\qed

	\begin{Remark}
		Notice that the proof of the first equality in Proposition \ref{prop:reformulation} does not depend on any special structure conditions on $\Om$ as in the framework of \cite{bouchard2016consistent}.
		In other words, it holds still true for an abstract space $(\Om, \Fc)$ with an arbitrary family of probability measures $\Pc$.
	\end{Remark}

	\begin{Remark}
		Let us observe that the reformulations in Proposition \ref{prop:reformulation} on the enlarged space do  not exactly correspond to  standard quasi-sure super-hedging problem. 
		Indeed, we still restrict the class of strategies to $\Fbb^{0}$-predictable processes, as opposed to $\Fbb$-predictable processes. 
		The fact that the formulation with these two different filtrations are equivalent will be proved by using a minimax argument in the next section.
	\end{Remark}

\subsubsection{Proof of Theorem \ref{thm:main_surrep}, case $e=0$}

	In view of Propositions \ref{prop:equiv_pricing} and \ref{prop:reformulation}, 
	 Theorem \ref{thm:main_surrep} will be proved if one can show that, with $g := \xi \cdot X_T$:
	$$
		\inf \Big\{
			y  \in \R~: y + (H \circ X)_T \ge g ~\Pcb_{\mathrm{int}} \mbox{-q.s.}, \mbox{ for some }H \in \Hc
		\Big\}
		~=~ 
		\sup_{\Qb \in \Qcb_0} \E^{\Qb} \big[ g \big].
	$$
	Let us start with a weak duality result,  which is an immediate consequence of   \cite[Lemmas A.2 and A.3]{BouchardNutz.13}.

	\begin{Lemma} \label{lemm:weak_hedging_duality}
		For any universally measurable variables $g: \Omb \to \R$, one has
		\b*
			&&
			\sup_{\Qb \in \Qcb_0} \E^{\Qb} \big[ g \big]
			~=~
			\sup_{\Qb \in \Qcb^{\mathrm{loc}}_0} \E^{\Qb} \big[ g \big] \\
			&\le&
			\inf \Big\{
				y  \in \R~: y + (H \circ X)_T \ge g ~\Pcb_{\mathrm{int}} \mbox{-q.s.}, \mbox{ for some }H \in \Hc
			\Big\}.
		\e*
	\end{Lemma}

	We   prove the converse inequality in the rest of this section.
	Let us proceed by induction, by first considering the  one period case $T=1$. 
	Recall that $\Lambda_0^{\mathrm{int}}(\om_0)$,
	$\Pcb^{\mathrm{int}, \delta}_0(\theta_0)$ and $\Qcb^{\delta}_0(\theta_{0})$
	are defined in Remark \ref{rem: Na loc condi}.

	\begin{Lemma} \label{lemm:surrepT1}
		Let $e=0$, $T=1$ and $g_1: \Omb \to \R \cup \{\infty\}$ be upper semi-analytic and such that
		$(\om, \theta_0, \theta_1)\in \Om\x \Lambda_{1}\x \Lambda_{1}\to g_1(\om, \theta_0, \theta_1)$ depends only on $(\om, \theta_1)$.
		Assume that $\NA(\Pcb_{\mathrm{int}})$ holds true.
		Then,
		\begin{align}
			\sup_{\theta_{0} \in \Lambda_{0}^{\rm int}(\om_{0})} \sup_{\Qb \in \Qcb^{\delta}_0(\theta_{0})} \E^{\Qb} \big[g_1 \big] &=\sup_{\Qb \in \Qcb_0} \E^{\Qb} \big[ g_1 \big]
			\label{eq:pric_hedg_duality_equiv}\\
			&=
			\inf \big\{ y \in \R ~: y + (H \circ X)_T \ge g_1, ~\Pcb_{\mathrm{int}} \mbox{-q.s.},~ H \in \Hc \big\}\nonumber\\
			&>-\infty.\nonumber
		\end{align}	
	\end{Lemma}
	\proof First, notice that $\Hc = \R^d$ when $T=1$,
	and for all $\theta_{0} \in \Lambda_{0}^{\rm int}(\om_{0})$,
	$$
		\big\{ \Pb \circ (g_1, X_1)^{-1} : \Pb \in \Pcb^{\mathrm{int},\delta}_{0}(\theta_{0}) \big\}
		~=~
		\big\{ \Pb \circ (g_1, X_1)^{-1} : \Pb \in \Pcb^{\mathrm{int},\delta}_{0}(\1) \big\},
	$$
	where $\1$ represents the vector of $\R^{d}$ with all entries equal to $1$.
	Then
	\begin{align} \label{eq:suprep_value}
		&
		\inf \big\{ y ~: y + (H \circ X)_1 \ge g_1, ~\Pcb_{\mathrm{int}} \mbox{-q.s.},~ H \in \Hc \big\} \nonumber \\
		&=
		\inf_{h_1 \in \R^{d}} \sup_{\theta_{0}\in \Lambda_{0}^{\rm int}(\om_{0})} \sup_{\Pb \lll \Pcb^{\mathrm{int},\delta}_{0}(\theta_{0})} 
		\E^{\Pb} [g_1 - h_1 \cdot (X_1 - S_{0}\theta_{0})] \nonumber \\
		&=
		\sup_{\theta_{0}\in \Lambda_{0}^{\rm int}(\om_{0})}  \inf_{h_1 \in \R^{d}} 
		 \sup_{\Pb \lll \Pcb^{\mathrm{int},\delta}_{0}(\theta_{0})} 
		\E^{\Pb} [g_1 - h_1 \cdot (X_1 -S_{0}\theta_{0})]\nonumber \\
		&=
		\sup_{\theta_{0} \in \Lambda_{0}^{\rm int}(\om_{0})} \inf \big\{ y ~: y + H_{1} \cdot (X_1-X_{0}) \ge g_1, ~\Pcb^{\mathrm{int},\delta}_{0}(\theta_{0}) \mbox{-q.s.},~ H \in \Hc \big\}.~~
	\end{align}
	In the above, the second equality follows by the minimax theorem since 
	$$
		(\theta_{0}, h_1) 
		\mapsto 
		\sup_{\Pb \lll\Pcb^{\mathrm{int},\delta}_{0}(\theta_{0})} \!\!\! \E^{\Pb} [ g_1- h_1 \cdot (X_1 - S_{0}\theta_{0})] 
		=
		\sup_{\Pb \lll\Pcb^{\mathrm{int},\delta}_{0}(\1)} \!\!\! \E^{\Pb} [ g_1- h_1 \cdot (X_1 - S_{0}\theta_{0})]
	$$
	is linear $\theta_{0}$ and convex in $h_1$, while the infimum over $h_{1}$ is concave and therefore lower semicontinuous (in particular, one can replace $\Lambda_{0}^{\rm int}(\om_{0})$ by its closure, that is a compact set, in all the above terms).
	Observe that $\Qcb^{\delta}_0(\theta_{0})$ 
	 is nonempty for every $\theta_{0} \in \Lambda_{0}^{\rm int}(\om_{0})$, recall Remark \ref{rem: Na loc condi} and Theorem \ref{thm:FundThmNT}. 
Then, by the duality result in \cite[Theorem 3.4]{BouchardNutz.13}, the right-hand side of \eqref{eq:suprep_value}   has a finite negative part and is equal to
	$$
		\sup_{\theta_{0} \in \Lambda_{0}^{\rm int}(\om_{0})} \sup_{\Qb \in \Qcb^{\delta}_0(\theta_{0})} \E^{\Qb} \big[g_1 \big]
		~=~
		\sup_{\Qb \in \Qcb_0} \E^{\Qb} \big[g_1 \big],
	$$
	where the last equality follows from the fact that $\Qc^{\delta}_0(\theta_{0}) \subset \Qc_0$ and that 
	every probability measure $\Qb$ in $\Qcb_0$ can be disintegrated into 
	a combinaison of elements in $(\Qc^{\delta}_0(\theta_{0}))_{\theta_{0} \in \Lambda_{0}^{\rm int}(\om_{0})}$.
	\qed

	\vspace{2mm}
	
	We now prepare for the general case $T\ge 1$, which is  based on a dynamic programming argument. 
	We extend the definitions of  $\Lambda_0^{\mathrm{int}}(\om_0)$, $\Pcb^{\mathrm{int}, \delta}_0(\theta_0)$ and $\Qcb^{\delta}_0(\theta_{0})$, see Remark \ref{rem: Na loc condi}, to an arbitrary initial time $t$ and initial path $\omb^{t}$. 
	For $t \ge 1$ and $\omb = \omb^t = (\om^t, \theta^t) \in \Omb_t$, we firt recall the definition of $\Lambda^{\rm int}_t (\om^{t}) $ that was already used in the proof of Proposition \ref{prop:equiv_NA}:
	$$
		\Lambda^{\rm int}_t (\om^{t}) 
		~:=~ 
		\{\theta_t \in \Lambda_1 ~: S_t(\om^t) \theta_t \in \mathrm{int}K^*_t(\om^t) \}
		~\subset~
		\Lambda_1.
	$$
	Next,
	recall that $\Pcb^{\mathrm{int}}_t(\omb) \subset \Bf(\Om_1 \x \Lambda_1)$ is  defined in \eqref{eq:def_Pcb_t} and note that it can also be written as $\Pcb^{\mathrm{int}}_t(\om)$ since it depends only on $\om$. We further define 
	$$\Pcb^{\mathrm{int},\delta}_t(\omb)  := \big\{
		\delta_{\omb^{t}} \otimes \Pb_{t+1} ~:  \Pb_{t+1} \in \Pcb^{\mathrm{int}}_t(\omb) \big\}$$
	and
	$$
		\Qcb^\delta_t(\omb)
		~:=~
		\big\{
			\Qb_{t+1} \lll \Pcb^{\mathrm{int},\delta}_t(\omb) ~: \E^{\Qb_{t+1}}[ X_{t+1} - X_t ] = 0
		\big\},
	$$
	as well as 	
	$$
		\Pct^{\mathrm{int},\delta}_t(\om) 
		:=
		\big\{
			(\delta_{\om^{t}}\times \mu(d \theta^t)) \otimes \Pb_{t+1}
			~:  \Pb_{t+1} \in \Pcb^{\mathrm{int}}_t(\om), ~\mu 
				\in \Bf \big(\Lambda^{\mathrm{int}}_0(\om_0) \x \cdots \x \Lambda^{\mathrm{int}}_t(\om^t) \big)
		\big\}.
	$$
	Let $g_{t+1}: \Omb_{t+1} \to \R \cup \{\infty\}$ be an upper semi-analytic functional 
	and be such that  $g_{t+1}(\om^{t+1}, \theta_0, \cdots, \theta_{t+1})$ depends only on $(\om^{t+1}, \theta_{t+1})$.
	We define 
	\be \label{eq:def_g_n}
		g_t(\omb^{t}) ~:=~ \sup_{\Qb \in \Qcb^\delta_t(\omb)} \E^{\Qb}[g_{t+1}],
	\ee
	which depends only on $(\om^t, \theta_t)$ by Remark \ref{rem:g_prime} below,
	and then define 
	$$
		g'_t(\om^{t}, h_{t}) 
		:=\!
		\sup_{\theta_t \in \Lambda^{\rm int}_t(\om^{t})} \big \{ g_t(\om^{t}, \theta_{t}) - h_t \cdot S_t(\om) \theta_{t} \big\},
		\;h_{t}\in \R^{d}.
	$$

	\begin{Remark} \label{rem:g_prime}
		Let $\bar \om=(\om,\theta)$ and $\bar \om'=(\om',\theta')$ be such that  $\om^{t} = (\om')^{t}$ and $\theta_t = \theta^{'}_t$.
		Then, it follows from the  definition of $\Pcb^{\mathrm{int},\delta}_t(\omb)$ and $\Qcb^\delta_t(\omb)$  that
		$$
			\{ \Qb \circ (g_{t+1}, X_t, X_{t+1})^{-1} ~:\Qb \in \Qcb^\delta_t(\omb) \}
			=
			\{ \Qb \circ (g_{t+1}, X_t, X_{t+1})^{-1} ~:\Qb \in \Qcb^\delta_t(\omb') \}.
		$$
		Hence, $g_t(\omb^{t})$ depends only on $(\om^{t}, \theta_t)$ for $\omb^{t} = (\om^{t}, \theta_0, \cdots, \theta_t)$.
	\end{Remark}
	
	\begin{Remark}\label{rem: NA Pintdelta omega} 
		{\rm (i)} For a fixed $\omb \in \Omb_t$, we define $\NA(\Pcb^{\mathrm{int},\delta}_t(\omb))$ by
		$$
			h\cdot(X_{t+1}-X_{t}) \ge 0~ \Pcb^{\mathrm{int},\delta}_t(\omb)\mbox{-q.s.}
			~~~\Longrightarrow~~~
			h\cdot(X_{t+1}-X_{t})= 0~ \Pcb^{\mathrm{int},\delta}_t(\omb)\mbox{-q.s.},
		$$
		for every  $h\in  \R^{d}$.
		Then, it follows from \cite[Theorem 4.5]{BouchardNutz.13} and Lemma \ref{lemma:panalytic} that 
		$\NA(\Pcb_{\mathrm{int}})$ implies that $\NA(\Pcb^{\mathrm{int},\delta}_t( \omb))$ holds for  all $ \omb \in \Omb$ outside a $\Pcb_{\mathrm{int}}$-polar set.
		
		\vspace{1mm}
		
		\noindent {\rm (ii)} Now, for a fixed $\om \in \Om_t$, let us define $\NA(\Pct^{\mathrm{int},\delta}_t(\om))$ by
		$$
			h(X_t) \cdot(X_{t+1}-X_{t}) \ge 0~ \Pct^{\mathrm{int},\delta}_t(\omb)\mbox{-q.s.}
			~\Longrightarrow~
			h(X_t) \cdot(X_{t+1}-X_{t})= 0~ \Pct^{\mathrm{int},\delta}_t(\omb)\mbox{-q.s.},
		$$
		for every universally measurable functions  $h: \R^{d} \to \R^d$.
		Then by applying Proposition \ref{prop:equiv_NA} with $\Pc(t,\om)$ in place of $\Pc$, one obtains that 
		$\NA2(t, \om)$ defined in \eqref{eq:NA2t} is equivalent to $\NA(\Pct^{\mathrm{int},\delta}_t(\om))$.
		
		\vspace{1mm}
		
		\noindent {\rm (iii)}  
		 It follows from {\rm (ii)}, Lemma \ref{lemm:NA2_polar} and Proposition \ref{prop:equiv_NA} that, whenever $\NA(\Pcb_{\mathrm{int}})$ holds,
		$\NA(\Pct^{\mathrm{int},\delta}_t(\om))$ holds for all $\om$ outside a $\Pcb_{\mathrm{int}}$-polar (or simply $\Pc$-polar set) set.
		Moreover, by similar arguments as in Remark \ref{rem: Na loc condi}, 
		$\NA(\Pct^{\mathrm{int},\delta}_t(\om))$ implies $\NA(\Pcb^{\mathrm{int},\delta}_t(\om,\theta))$ for all $\theta\in \Lambda$.

	\end{Remark}

	\begin{Lemma} \label{lemm:selec_h_n1} Assume that  $\NA(\Pcb_{\mathrm{int}})$ holds. 
		Then, both $g'_t$ and $g_t$ are upper semi-analytic. 
		Moreover, 
		there is a universally measurable map $h_{t+1}: \Om_t \x \R^d \to \R^d$ and a $\Pc$-polar set $N$ such that, 
		for every $(\om, h_t) \in N^{c} \x \R^d$,
		one has $g'_t(\om^{t}, h_t) > -\infty$ and 		
		\b*
			g'_t(\om^{t}, h_t) + h_t \cdot X_t   + h_{t+1}(\om^{t}, h_t) \cdot (X_{t+1} - X_t ) 
			~\ge~
			g_{t+1}~~ \Pct^{\mathrm{int},\delta}_t(\om) \mbox{-q.s.}
		\e*
	\end{Lemma}
	\proof The proof follows from the same measurable selection arguments as in   \cite[Lemma 4.10]{BouchardNutz.13}.
	We provide a sketch of proof for completeness.
	 Let   define 
	 	\begin{equation*}
		\pi'_t (\om^{t}, h_t) 
		 := 
		\inf\big\{
		y \in \R  :   
			y + h_t \cdot X_t + h_{t+1}\cdot (X_{t+1} - X_t) 
			\ge
			g_{t+1}~
		 		\Pct^{\mathrm{int},\delta}_t(\om)\mbox{-q.s.},\; h_{t+1}\in \R^{d}
		\big\}.
	\end{equation*}
	Then, the  same minimax theorem argument as the one used for the proof of  Lemma \ref{lemm:surrepT1} implies that
	\b*
		\pi'_t(\om^{t}, h_t) = g'_t( \om^{t}, h_t)
		~>~ -\infty
		~~~\mbox{if}~\NA(\Pct^{\mathrm{int},\delta}_t(\om)) ~\mbox{holds true}. 
	\e*
	In view of {\rm (iii)} in Remark \ref{rem: NA Pintdelta omega}, this is true outside a $\Pc$-polar set $N$.
	
	Further, $g_{t+1}$ is assumed to be upper semi-analytic, $X_{t}$ is Borel measurable,
	the graph of $\Qcb^\delta_t$ is analytic by \cite[Lemma 4.8]{BouchardNutz.13},
	and the graph of $\mathrm{int}K^*_t$ is Borel.
	It thus follows from a   measurable selection argument (see e.g.~\cite[Propositions 7.26, 7.48]{BertsekasShreve.78}) that
	the maps $\omb^{t} \mapsto g_t(\omb^{t})$ and $(\om^{t}, h_t) \mapsto g'_t(\om^{t}, h_t)$ are
	both upper semi-analytic.

	We now define $\hat{g}_t:=g'_t \mathbf{1}_{\R}(g'_t)$,
	which is universally measurable, i.e.~in $\Uc(\Om_t \times \R^d)$, and consider the  random set: 
	$$
		\Psi(\om^{t}, h_t):=
		\big\{ y \in \R^d: \gamma_y(\cdot , h_t) \le 0		~\Pct^{\mathrm{int},\delta}_t(\om)
		\mbox{-q.s.} 
		\big\}
	$$
	where 
	\begin{align*}
	\gamma_y(\cdot , h_t) 
	:=&
	g_{t+1}-\hat{g}_t(\cdot, h_t)
	-y\cdot\big(X_{t+1} -X_t\big)  - h_t X_t.
\end{align*}
	It is enough to show that $\{ \Psi \neq \emptyset \}$ is universally measurable and that $\Psi$ admits a universally measurable selector $h_{t+1}(\cdot)$ on $\{ \Psi \neq \emptyset \}$. 
It is not hard to see that\footnote{$\mbox{USA}[\Uc(\Om_t \times \R^d) \otimes  \Bc(\R^d \times \Omb_1) ]$ denotes the convex cone generated by both the upper semianalytic maps on $(\Om_t \times \R^d)\x  (\R^d \times \Omb_1) $ and the  $\Uc(\Om_t \times \R^d) \otimes  \Bc(\R^d \times \Omb_1)$-measurable functions.} $\gamma_y \in \mbox{USA}[\Uc(\Om_t \times \R^d)\otimes\Bc(\R^d \times \Omb_1)]$. Note that, given a probability measure $\hat \P$, $\gamma_y(\cdot, h_t) \leq 0~~\hat{\P}$-a.s.~iff $\E^{\Pt}[\gamma_y(\cdot, h_t)] \leq 0 ~ \mbox{for~all~} \Pt \ll \hat{\P}$. By an application of  \cite[Lemma 3.2]{BouchardNutz.13}, it is further equivalent to have the above for all $\Pt \ll \hat{\P}$ satisfying $\E^{\Pt}[|X_{t+1}-X_t|] < \infty$ and $\E^{\Pt}[|X_t|]< \infty$. Therefore, we introduce the random set
$$
\Pct^{\mathrm{int},\delta'}_t(\om):=\{ \Pt \lll \Pct^{\mathrm{int},\delta}_t(\om) ~: \E^{\Pt}[|X_{t+1}-X_t|]+\E^{\Pt}[|X_t|]< \infty \}.
$$
By the same arguments in the proof of   \cite[Lemma 4.8]{BouchardNutz.13}, one can prove that $\Pct^{\mathrm{int},\delta'}_t$ has an analytic graph. Define now 
$$
\Gamma_y(\om^t, h_t):=\sup_{\Pt \in \Pct^{\mathrm{int},\delta'}_t(\om^t)} \E^{\Pt}[\gamma_y (\cdot, h_t)],
$$
so that $\gamma_y (\cdot, h_t) \leq 0~\Pct^{\mathrm{int},\delta}_t(\om)\mbox{-q.s.}$  iff $\Gamma_y(\om^t, h_t) \leq 0$. We  now show that $(\om^t, h_t) \mapsto \Gamma_y(\om^t, h_t)$ is universally measurable. Indeed, the first term in the difference
\begin{align*}
&\E^{\Pt}[g_{t+1}(\om^{t},\cdot)]-\hat{g}_t(\omega^{t},h_{t})-y\cdot \E^{\Pt}[X_{t+1}(\om^{t},\cdot)-X_t(\om^{t},\cdot)]-h_t\cdot  \E^{\Pt}[X_t(\om^{t},\cdot)] 
\end{align*}
is a upper semianalytic function of $(\om^{t},  \Pt)$. The second term is universally measurable. The  third and the fourth terms are Borel.  As a result, $(\om^{t}, h_t, \Pt) \mapsto \E^{\Pt}[\gamma_y(\om^{t},\cdot, h_t)]$ is in $\mbox{USA}[\Uc(\Om_t \times \R^d) \otimes \Bc(\Bf(\R^d \times \Omb_1))]$. Thus, by the Projection Theorem in the form of  \cite[Lemma 4.11]{BouchardNutz.13}, 
\begin{align*}
\{ \Gamma_y > c \} &= \mbox{proj}_{\Om_t \times \R^d} \{ (\om^t, h_t, \Pt) : (\om^t,\Pt) \in \mbox{graph}(\Pct^{\mathrm{int},\delta'}_t),~ \E^{\Pt}[\gamma_y(\om^{t},\cdot, h_t)] > c \} \\
&\in \Uc(\Om_t \times \R^d)
\end{align*}
for all $c \in \R$. This means that  $(\om^{t}, h_t) \mapsto \Gamma_y(\om^{t}, h_t)$ is universally measurable for any fixed $y$.

On the other hand, given  $(\om^{t}, h_t) \in \Om_t \times \R^d$  and $m \geq 1$,  the function $y \mapsto \Gamma_y(\om^{t}, h_t) \wedge m$ is lower semicontinuous as  the supremum over   $\Pct^{\mathrm{int},\delta'}_t$ of a family of    continuous functions. By   \cite[Lemma 4.12]{BouchardNutz.13},  $(\om^{t}, h_t, y) \mapsto \Gamma_y(\om^{t}, h_t) \wedge m$ is $\Uc(\Om_t \times \R^d) \otimes \Bc(\R^d)$-measurable as well. As a result,   
\b*
	\mbox{graph}(\Psi)
	~=~
	\{ (\om^{t}, h_t, y): \Gamma_y(\om^{t}, h_t) \leq 0 \} 
	&\in&
	\Uc(\Om_t \times \R^d) \otimes \Bc(\R^d) \\
	&&~~~
	\subset
	\mbox{A}[\Uc(\Om_t \times \R^d) \otimes \Bc(\R^d)],
\e*
where $ \mbox{A}[\Uc(\Om_t \times \R^d) \otimes \Bc(\R^d)]$ stands for the collection analytic sets of $\Uc(\Om_t \times \R^d) \otimes \Bc(\R^d)$. 
Finally, \cite[Lemma 4.11]{BouchardNutz.13} yields that $\{ \Psi \neq \emptyset \} \in \Uc(\Om_t \times \R^d)$ and that $\Psi$ admits a $\Uc(\Om_t \times \R^d)$-measurable selector on $\{ \Psi \neq \emptyset \} $.
	\qed
	
	\vspace{2mm}

	{\bf \noindent Proof of Theorem \ref{thm:main_surrep} (case $e=0$).}
	The existence of the optimal super-hedging strategy will be proved in Lemma \ref{lemm:exist_hedging} below for the  general case $e\ge 1$.
	By Propositions \ref{prop:equiv_pricing} and \ref{prop:reformulation}  and Lemma \ref{lemm:weak_hedging_duality},
	it is enough to prove
	\be\label{eq:pric_hedg_duality_equiv_loc}
		\sup_{\Qb \in \Qcb^{\mathrm{loc}}_0} \E^{\Qb} \big[ g \big]
		\ge
		\inf \Big\{
			y  \in \R~: y + (H \circ X)_T \ge g ~\Pcb_{\mathrm{int}} \mbox{-q.s.}, \mbox{ for some }H \in \Hc
		\Big\},~~
	\ee 
	for $g := \xi \cdot X_T$.
	We use an induction argument. 
	Recall  that \eqref{eq:pric_hedg_duality_equiv_loc} is already proved for the case   $T=1$, this is the content of Lemma \ref{lemm:surrepT1}. 
	Assume that \eqref{eq:pric_hedg_duality_equiv_loc} holds true for  $T = t$ and let us   prove that it also holds true for the case $T=t+1$. 

	Given an upper semianalytic random variable $g_{t+1} := \Omb_{t+1} \to \R \cup \{ \infty\}$ such that
	$g_{t+1}(\om^{t+1}, \theta_0, \cdots, \theta_{t+1})$ depends only on $(\om^{t+1}, \theta_{t+1})$.
	we define $g_t$ by \eqref{eq:def_g_n}, and denote
	$$
		\pi_0^t (g_t)
		~:=~
		\inf \big\{ y ~: y + (H \circ X)_t \ge g_t~~\Pcb_{\mathrm{int}} \mbox{-q.s.},~ H \in \Hc \big\}.
	$$
	Fix $(y, H) \in \R \x \Hc$ 
	such that $y + (H \circ X)_t \ge g_t, ~\Pcb_{\mathrm{int}} \mbox{-q.s.}$. 
	Then,  $y + (H \circ X)_{t-1} - H_{t}\cdot X_{t-1} \ge g_t - H_{t}\cdot X_{t}$ $\Pcb_{\mathrm{int}} \mbox{-q.s.}$
	and therefore 
	$$
		y + (H \circ X)_{t-1} - H_{t}\cdot X_{  t-1}\ge g'_{t}(\cdot, H_{t})~\Pcb_{\mathrm{int}} \mbox{-q.s.}
	$$
	Hence, if we 
	define  $H'$ by $H'_s := H_s$ for $s \le t$ and	
	$H'_{t+1}(\om^{t}) := h_{t+1}(\om^{t}, H_t(\om^{t-1}))$, with $h_{t+1}$ as in Lemma \ref{lemm:selec_h_n1}, we  obtain 
	$$
		y + (H' \circ X)_{t+1} \ge g_{t+1}~~\Pcb_{\mathrm{int}} \mbox{-q.s.}
	$$
	Hence,
	$$
		\pi_0^{t+1} (g_{t+1}) 
		~\le~ 
		\pi_0^t (g_t)
		~=~
		\sup_{\Qb \in \Qcb_0^{\mathrm{loc}}} \E^{\Qb}\big[ g_t \big]
		~\le~
		\sup_{\Qb \in \Qcb_0^{\mathrm{loc}}} \E^{\Qb}\big[ g_{t+1} \big],	
	$$
	where the last inequality follows from a classical concatenation argument.
	This is in fact \eqref{eq:pric_hedg_duality_equiv_loc} for the case $T= t+1$, and we hence conclude the proof.
	\qed

\subsection{Proof of Theorem \ref{thm:main_surrep}: case $e \ge 1$}

	To take into account the transaction costs generated by the trading of the static options $(\zeta_i, i=1, \cdots, e)$, we introduce a further enlarged space:
	$$
		\widehat \Lambda := \prod_{i=1}^e [-c_i, c_i],
		~~
		\Omh := \Omb \x \widehat \Lambda ,
		~~
		\Fch_t := \Fcb_t \otimes \Bc\big(\widehat \Lambda \big),
		~~
		\Pch_{\mathrm{int}} := \big\{
		\Ph \in \Bf(\Omh)~: \Ph|_{\Omb} \in \Pcb_{\mathrm{int}}
		\big\},
	$$
	and define
	$$
		\hat f_i : \Omh ~\longrightarrow~ \R ,
		~~~
		\hat f_i(\omh) = \zeta_{i}(\om) \cdot X_T(\omb) - \hat \theta_i
		~~\mbox{for all}~\omh = (\omb, \hat \theta) = (\om, \theta, \hat \theta) \in \Omh.
	$$
	The process $(X_t)_{0 \le t \le T}$ and the random variable $g := \xi \cdot X_T$ defined on $\Omb$ can be naturally extended on $\Omh$.
	We can then consider the super-hedging problem on $\Omh$:
	$$
		\hat \pi_e(g) 
		~:=~
		\inf\Big\{ 
			y ~: y + \sum_{i=1}^e \ell_i \hat f_i + (H \circ X)_T \ge g, ~\Pch_{\mathrm{int}} \mbox{-q.s.}, ~\ell \in \R^e,~ H \in \Hc
		\Big\}.
	$$
	Let us also introduce
	$$
		\Qch_e 
		~:=~ 
		\big\{ 
			\Qh \in \Bf(\Omh) ~: \Qh \lll \Pch_{\mathrm{int}}, ~ X ~\mbox{is}~(\Fbh, \Qh) \mbox{-martingale}, ~\E^{\Qh}[\hat f_i] = 0,~ i=1,\cdots, e
		\big\},
	$$
	and
	$$
		\Qch_e^{\varphi} := \big\{
			\Qh \in \Qch_e ~: \E^{\Qh}[\varphi] < \infty
		\big\},
		~~~\mbox{for}~~
		\varphi := 1 + |g| + \sum_{i=1}^e | \hat f_i|.
	$$

	\begin{Lemma} \label{lemm:exist_hedging} Let  $\NA2(\Pc)$ hold. Assume further that and \eqref{eq: NA avec option et CT} holds true for all $\ell \in \R^{e}$ and $\eta\in \Ac$. Then:

		\noindent (a) There exist $\hat \ell \in \R^e$ and a $\F$-predictable process $\widehat H$ such that
		\be \label{eq:exist_hedging}
			\hat \pi_e(g) + \sum_{i=1}^e \hat \ell_i \hat f_i + (\widehat H \circ X)_T \ge g, ~\Pch_{\mathrm{int}} \mbox{-q.s.}
		\ee
		(b) Consequently, there exists  $(\hat \eta, \hat \ell) \in \Ac \x \R^e$ such that 
		\be \label{eq:exist_hedging avec CT}
			\pi_e(\xi){\rm 1}_{d} 
			+ \sum_{i=1}^e  \Big(  \hat \ell_i \zeta_i -  |\hat \ell_i| c_i {\rm 1}_{d} \Big) 
			+ \sum_{t =0}^T \hat \eta_t - \xi \in K_T,
			~\Pc \mbox{-q.s.}
		\ee
	\end{Lemma}
	\proof (a) It suffices to show that the collection of claims that can be super-hedged from $0$ is closed for the $\widehat \Pc_{\mathrm{int}}$-q.s.~convergence. Note the results in  \cite[Section 2]{BouchardNutz.13} are given in a general abstract context, where the underlying asset is not assumed to be adapted to the filtration of the strategy. 
	Then,  \cite[Theorem 2.3]{BouchardNutz.13} implies our claim in the case $e=0$ (recall that  $\NA2(\Pc)$ implies $\NA(\Pcb_{\rm int})$). Assume now that it holds for $e-1\ge 0$ and let us deduce that it holds for $e$ as well.

	Let $(g^{n})_{n\ge 1}\subset \L^{0}$ be such that $g^{n}\to g$  $\widehat \Pc_{\mathrm{int}}$-q.s., and $\pi_{e}(g^{n})\le 0$ for $n\ge 1$. Let $(\hat \ell^{n})_{n} \subset  \R^{e}$ and let $(\widehat H^{n})_{n\ge 1}$ be a sequence of  $\F$-predictable processes such that 
	$  \sum_{i=1}^e \hat \ell^{n}_i \hat f_i + (\widehat H^{n} \circ X)_T \ge g^{n}$  $\Pch_{\mathrm{int}}\mbox{-q.s.}$ 
	If $(\hat \ell^{n}_{e})_{n\ge 1}$ is bounded, then one can assume that it converges to some $\hat \ell_{e}\in \R$. 
	Hence, \cite[Theorem 2.3]{BouchardNutz.13} implies that one can find $\hat \ell \in \R^{e-1}$ and a  $\F$-predictable process $\widehat H$ such that $ \sum_{i=1}^{e-1}   \hat \ell_i \hat f_i + ( \widehat H \circ X)_T \ge g - \hat \ell_{e}\hat f_{e}$ $\Pch_{\mathrm{int}} \mbox{-q.s.}$
	
	If $(\hat \ell^{n}_{e})_{n\ge 1}$ is not bounded, then one can assume that $|\hat \ell^{n}_{e}|\to \infty$, 
	so that $(g^{n}- \hat \ell^{n}_{e} \hat f_{e})/(1+|\hat \ell^{n}_{e}|)$ $\to$ $-\chi \hat f_{e}$ $\widehat \Pc_{\mathrm{int}}$-q.s.~for some $\chi \in \{-1,1\}$ 
	and \cite[Theorem 2.3]{BouchardNutz.13} implies that one can find $\hat \ell\in \R^{e-1}$ and a  $\F$-predictable process $\widehat H$ such that $\sum_{i=1}^{e-1}  \hat \ell_i \hat f_i + ( \widehat H \circ X)_T \ge -\chi \hat f_{e}$ $\Pch_{\mathrm{int}}\mbox{-q.s.}$  
	By similar arguments as in the proof of Proposition \ref{prop:reformulation}, this implies that  
	$\chi \hat f_{e}-|\chi|c_{e}{\rm 1}_{d}+ \sum_{i=1}^{e-1}   (\hat \ell_i  \zeta_i - |\hat \ell_{i}|c_{i}{\rm 1}_{d})+\sum_{t=0}^{T} \eta_{t}\in K_{T}~\Pc\mbox{-q.s.}$ for some $\eta \in \Ac$. Then, $\chi=0$ by \eqref{eq: NA avec option et CT}, a contradiction.
	
	\vspace{2mm}
	
	\noindent (b)
	Finally, by the same arguments as in the proof of Proposition \ref{prop:reformulation},
	one can show $\pi_e(\xi) = \hat \pi_e(g)$ for $g := \xi \cdot X_T$.
	Moreover, using the construction \eqref{eq:def_eta_d},
	one can obtain explicitly $(\hat \eta, \hat \ell)$ satisfying \eqref{eq:exist_hedging avec CT}
	from $(\widehat H, \hat \ell)$ satisfying \eqref{eq:exist_hedging}.
	\qed
	
	\vspace{2mm}

	{\bf \noindent Proof of Theorem \ref{thm:main_surrep} (case $e \ge 1$).}
	The existence of a  super-hedging strategy has been  proved in Lemma \ref{lemm:exist_hedging} above.
	Moreover, it is easy to adapt the arguments of Propositions \ref{prop:equiv_pricing} and \ref{prop:reformulation} to obtain 
	$$
	  	\pi_e(\xi)=\hat \pi_e(g)
		\;\mbox{ and } \;
	   	\sup_{(\Q,Z) \in \CPSt_e} \E^{\Q} \big[\xi \cdot Z_T \big]
		=
		\sup_{\Qh \in \Qch_e} \E^{\Qh} \big[ g \big],
		~~~\mbox{for}~g:=\xi \cdot X_T.
	$$
	Remember that, by \cite[Lemma A.3]{BouchardNutz.13}, one has $\sup_{\Qh \in \Qch^{\varphi}_e} \E^{\Qh}[ g] 
	= \sup_{\Qh \in \Qch_e} \E^{\Qh}[ g]$.
	Hence,  it is enough to prove that
	\be \label{eq:dualityOmh}
		\hat \pi_e(g)
		~=~
		\sup_{\Qh \in \Qch^{\varphi}_e} \E^{\Qh}[ g].
	\ee
	Note that we have already proved \eqref{eq:dualityOmh} for the case $e=0$ in Section \ref{subsec:proof_e0}, although the formulations are slightly different (the additional randomness induced by $\widehat \Lambda$ obviously does not play any role when $e=0$).
	We argue by induction as in the proof of  \cite[ Theorem 5.1]{BouchardNutz.13}. 
	Let us assume that \eqref{eq:dualityOmh} holds for $e-1\ge 0$ and then prove it  for $e$.
	
	First, it follows from \eqref{eq: NA avec option et CT} that we  can not find $\eta \in \Ac$ and $(\ell_i)_{1 \le i \le e} \neq 0$ such that $\sum_{i=1}^{e}   (\ell_i  \zeta_i-|\ell_{i}|c_{i}{\rm 1}_{d})+\sum_{t=0}^{T} \eta_{t}\in K_{T}$ $\Pc\mbox{-q.s.}$
	By the same arguments as in the proof of Proposition \ref{prop:reformulation},
	there is no $H \in \Hc$, $\ell_1, \cdots, \ell_{e-1}$ and $\ell_{e}\in \{-1,1\}$ such that $ \sum_{i=1}^{e-1} \ell_i \hat f_i  + (H \circ X)_T\ge -\ell_{e}\hat f_{e}$, $\Pch_{\mathrm{int}}$-q.s.
	It follows that $\hat \pi_{e-1}(\hat f_{e}),\hat \pi_{e-1}(-\hat f_{e})>0$, which, by Lemma \ref{lemm:exist_hedging} and our induction hypothesis, implies that 
	there is $\Qh_{-}, \Qh_+ \in \Qch^{\varphi}_{e-1}$ such that
	\be\label{eq: pie de fe+1 avec 0 au milieu}
		- \hat \pi_{e-1}(- \hat f_{e}) ~<~ \E^{\Qh_-}[\hat  f_{e}] 
		~<~ 0 ~<~  \E^{\Qh_+}[ \hat f_{e}] ~<~ \hat \pi_{e-1}(\hat f_{e}).
	\ee
	We now  claim that 
	\be \label{eq:claim_approx_sup_duality}
		\mbox{there exists a sequence}~ \big(\Qh_n \big)_{n\ge 1} \subset \Qch^{\varphi}_{e-1} 
		~\mbox{s.t.}~ \E^{\Qh_n}[\hat f_{e}] \to 0, ~\E^{\Qh_n}[ g] \to \hat  \pi_{e}(g).~~~
	\ee
	Indeed, if the above fails, then 
	$$
		(0, \hat \pi_{e}(g) ) \notin \overline{\big\{ \E^{\Qh}[(\hat f_{e}, g)] : \Qh \in \Qch^{\varphi}_{e-1} \big\}} \subset \R^2,
	$$
	and one obtains a contradiction by following  line by line the same arguments in the end of the proof  of \cite[Theorem 5.1]{BouchardNutz.13}.
	In view of \eqref{eq: pie de fe+1 avec 0 au milieu} and \eqref{eq:claim_approx_sup_duality}, we can find $(\lambda_n^-, \lambda_n, \lambda_n^+) \in [0,1]$ such that $\lambda_n^- + \lambda_n +\lambda_n^+ = 1$, $(\lambda_n^-, \lambda_n^+) \to 0$, and
	$$
		\Qh'_n := \lambda_n^- \Qh_- + \lambda_n \Qh_n + \lambda_n^+ \Qh_+ \in \Qch^{\varphi}_{e-1}
		~~\mbox{satisfies}~
		\E^{\Qh'_n} [f_{e}] = 0.
	$$
	In particular, one has $\Qh'_n \in \Qch^{\varphi}_e$ and hence $\Qch^{\varphi}_e$ is nonempty,
	which implies that $\CPSt_e$ is nonempty by the projection argument in Proposition \ref{prop:equiv_pricing}.
	
	Moreover, since $\E^{\Qh'_n} [ g] \to \hat \pi_{e}(g)$, this shows that 
	$$
		\sup_{\Qh \in \Qch^{\varphi}_{e}} \E^{\Qh}[ g] 
		~\ge~
		\hat \pi_{e}(g).
	$$
	To conclude, it is enough to notice that the reverse inequality is the classical weak duality which can be easily obtained from \cite[Lemmas A.1 and A.2]{BouchardNutz.13}.
	\qed

\end{document}